\documentclass[a4paper,twoside]{amsart}

\usepackage{amsmath}
\usepackage{amsthm}
\usepackage{amsfonts}
\usepackage{amssymb}

\allowdisplaybreaks  

\newtheorem*{gp}{The geometry problem}
\newtheorem*{theo}{Theorem}

\newtheorem{theorem}{Theorem}[section]
\newtheorem{lemma}[theorem]{Lemma}
\newtheorem{corollary}[theorem]{Corollary}
\newtheorem{proposition}[theorem]{Proposition}

\theoremstyle{definition}
\newtheorem{definition}[theorem]{Definition}
\newtheorem{example}[theorem]{Example}
\newtheorem{remark}[theorem]{Remark}

\DeclareMathOperator{\Ad}{Ad}
\DeclareMathOperator{\ad}{ad}
\DeclareMathOperator{\Hom}{Hom}

\newcommand{\be}{\begin{equation}}
\newcommand{\ee}{\end{equation}}
\newcommand{\ben}{\begin{equation*}}
\newcommand{\een}{\end{equation*}}
\newcommand{\bal}{\begin{aligned}}
\newcommand{\eal}{\end{aligned}}
\newcommand{\bma}{\begin{pmatrix}}
\newcommand{\ema}{\end{pmatrix}}

\newcommand{\inv}[2]{\mathbf{#1}_{\mathbf{#2}}}

\newcommand{\wt}[1]{\widetilde{#1}}
\newcommand{\wh}[1]{\widehat{#1}}
\newcommand{\goth}[1]{\mathfrak{#1}}

\newcommand{\inc}[2]{\mathbf{#1}^c_{\mathbf{#2}}}
\newcommand{\hc}[1]{\theta^{#1}}
\newcommand{\vc}[1]{\Omega_{#1}}

\newcommand{\inp}[2]{\mathbf{#1}^p_{\mathbf{#2}}}
\newcommand{\hp}[1]{\theta^{#1}}
\newcommand{\vp}[1]{\Omega_{#1}}

\renewcommand{\inf}[2]{\mathbf{#1}^f_{\mathbf{#2}}}
\newcommand{\hf}[1]{\theta^{#1}}
\newcommand{\vf}[1]{\Omega_{#1}}

\newcommand{\der}{{\rm d}}
\newcommand{\w}{{\scriptstyle\wedge}\,}
\newcommand{\hook}{\lrcorner}
\newcommand{\semi}[1]{\oplus_{#1}}

\newcommand{\C}{\mathcal{C}}
\newcommand{\co}{\goth{co}}
\newcommand{\conf}{\goth{conf}}
\newcommand{\D}{\mathcal{D}}
\newcommand{\Der}{\goth{D}}

\newcommand{\g}{\goth{g}}
\renewcommand{\gg}{\mathbf{g}}
\newcommand{\gl}{\goth{gl}}
\renewcommand{\H}{\mathcal{H}}
\newcommand{\h}{\goth{h}}
\newcommand{\J}{{J}}

\newcommand{\M}{{M}}

\newcommand{\m}{\goth{m}}
\renewcommand{\O}{\mathcal{O}}
\renewcommand{\o}{\goth{o}}
\renewcommand{\P}{{P}}

\newcommand{\R}{{\rm R}}
\newcommand{\real}{\mathbb{R}}
\newcommand{\Ric}{{\rm Ric}}
\renewcommand{\S}{\mathcal{S}}

\renewcommand{\sp}{\goth{sp}}

\renewcommand{\u}{\goth{u}}
\newcommand{\V}{{\mathcal V}}

\begin{document}

\title{Geometry of third-order ODEs}

\author{Michal Godlinski}
\address{Instytut Matematyczny Polskiej Akademii Nauk, ul. Sniadeckich 8, Warszawa, Poland}
\email{godlinski@impan.gov.pl}

\author{Pawel Nurowski}
\address{Instytut Fizyki Teoretycznej,
Uniwersytet Warszawski, ul. Hoza 69, Warszawa, Poland}
\email{nurowski@fuw.edu.pl}
\date{24 February 2009}

\begin{abstract}
We address the problem of local geometry of third order ODEs
modulo contact, point and fibre-preserving transformations of
variables. Several new and already known geometries are described
in a uniform manner by the Cartan method of equivalence. This
includes conformal, Weyl and metric geometries in three and six
dimensions and contact projective geometry in dimension three.
Respective connections for these geometries are given and their
curvatures are expressed by contact, point or fibre-preserving
relative invariants of the ODEs.
\end{abstract}

\maketitle

\tableofcontents

\section{Introduction}\label{ch.problem}
 \noindent This paper addresses the problem of geometry
of third order ordinary differential equations (ODEs) which is
stated as follows.
\begin{gp}
Determine geometric structures defined by a class of equations
\ben y'''=F(x,y,y',y'') \een equivalent under certain type of
transformations. Find relations between invariants of the ODEs and
invariants of the geometric structures.
\end{gp}

One may consider equivalence with respect to several types of
transformations, in this work we focus on three best known types:
contact, point and fibre-preserving transformations. The
fibre-preserving transformations are those which transform the
independent variable $x$ and the dependent variable $y$ in such a
way that the notion of the independent variable is retained, that
is the transformation of $x$ is a function of $x$ only:
\be\label{e.fp}
   x\mapsto\bar{x}=\chi(x),\quad\quad\quad\quad  y\mapsto\bar{y}=\phi(x,y).
\ee The transformation rules for the derivatives are already
uniquely defined by above formulae. Let us define the total
derivative to be \ben
  \Der=\partial_x+y'\partial_y+y''\partial_{y'}+y'''\partial_{y''}.
\een Then
\begin{subequations}\label{e.tprol}
\begin{align}
 y'&\mapsto\frac{\der\bar{y}}{\der\bar{x}}=\frac{\Der \phi}{\Der \chi}, \label{e.tprol1}  \\
 y''&\mapsto\frac{\der^2\bar{y}}{\der\bar{x}^2}= \frac{\Der}{\Der\chi}\left(\frac{\Der\phi}{\Der \chi}\right),\label{e.tprol2} \\
 y'''&\mapsto \frac{\der^3\bar{y}}{\der\bar{x}^3}=\frac{\Der}{\Der\chi}\left(\frac{\Der}{\Der\chi}\left(\frac{\Der\phi}{\Der \chi}\right)\right).\label{e.tprol3}
\end{align}
\end{subequations}
The point transformations of variables mix $x$ and $y$ in an
arbitrary way \be \label{e.point}
   x\mapsto\bar{x}=\chi(x,y),\quad\quad\quad\quad  y\mapsto\bar{y}=\phi(x,y),
\ee with the derivatives transforming as in \eqref{e.tprol}. The
contact transformations are more general yet. Not only they
augment the independent and the dependent variables but also the
first derivative \be\label{e.cont}
\begin{aligned}
   x&\mapsto\bar{x}=\chi(x,y,y'),\notag\\
   y&\mapsto\bar{y}=\phi(x,y,y'),\\
   y'&\mapsto\frac{\der \bar{y}}{\der\bar{x}}=\psi(x,y,y').\notag
\end{aligned}
\ee However, the functions $\chi$, $\phi$ and $\psi$ are not
arbitrary here but subjecting to \eqref{e.tprol1} which now yields
two additional constraints \ben
\psi=\frac{\Der\phi}{\Der\chi}\quad \iff \quad
\bal &\psi\chi_{y'}=\phi_{y'}, \\
  &\psi(\chi_x+y'\chi_y)=\phi_x+y'\phi_y, \eal
\een guaranteeing that $\der\bar{y}/\der\bar{x}$ really transforms
like first derivative. With these conditions fulfilled second and
third derivative transform through \eqref{e.tprol2} --
\eqref{e.tprol3}. We always assume in this work that ODEs are
defined locally by a smooth real function $F$ and are considered
apart from singularities. The transformations are always assumed
to be local
diffeomorphisms. %Particularly important in study of the geometry
%problem are relative invariants defined as follows.

A pioneering work on geometry of ODEs of arbitrary order is Karl
W\"unschmann's PhD thesis \cite{Wun} written under supervision of
F. Engel in 1905. In this paper K. W\"unschmann observed that
solutions of an $n$th-order ODE
$y^{(n)}=F(x,y,y',\ldots,y^{(n-1)})$ may be considered as both
curves $y=y(x,c_0,c_1,\ldots,c_{n-1})$ in the $xy$ space and
points $c=(c_0,\ldots,c_{n-1})$ in the solution space $\real^n$
parameterized by values of the constants of integration $c_i$. He
defined a relation of $k$th-order contact between infinitesimally
close solutions considered as curves; two solutions $y(x)$ and
$y(x)+\der y(x)$ corresponding to $c$ and $c+\der c$ have the
$k$th-order contact if their $k$th jets coincide at some point
$(x_0, y_0)$. W\"unschmann's main question was how the property of
having $(n-2)$nd contact for $n=3,4$ and $5$ might be described in
terms of the solution space. In particular he examined third-order
ODEs and showed that there is a distinguished class of ODEs
satisfying certain condition for the function $F$, which we call
the W\"unschmann condition. For a third-order ODE in this class,
the condition of having first order contact is described by a
second order Monge equation for $\der c$.  This Monge equation is
nothing but the condition that the vector defined by two
infinitesimally close points $c$ and $c+\der c$ is null with
respect to a Lorentzian conformal metric on the solution space.
The last observation, although not contained in W\"unschmann's
work, follows immediately from his reasoning and was later made by
S.-S. Chern \cite{Chern}, who cited W\"unschmann's thesis.

The main contribution to the issue of point and contact geometry
of third-order ODEs was made by E. Cartan and S.-S. Chern in their
classical papers \cite{Car1, Car2, Car3} and \cite{Chern}. E.
Cartan \cite{Car2} considered a third-order ODE modulo point
transformations and, applying his method of equivalence,
constructed a 7-dimensional manifold $\P$ together with a fixed
coframe
$\theta^1,\theta^2,\theta^3,\theta^4,\Omega_1,\Omega_2,\Omega_3$
which encodes all the point invariant information about the ODE.
This means that two ODEs are point equivalent if and only if their
associated coframes are diffeomorphic. In the language of
contemporary differential geometry $\P$ is a principal bundle
$H_3\to\P\to\J^3$ over the second jet space with the structure
group $H_3=\real\times(\real\ltimes\real)$ while the coframe
defines a $\co(2,1)\semi{.}\real^3$-valued\footnote{Notation is
explained on pages \pageref{s.notation} -- \pageref{ch.prel}.}
Cartan connection on $\P$. The point invariant information about
the ODE is contained in the curvature of this connection and its
coframe derivatives. E. Cartan observed that some equations have a
nontrivial geometry on their solution spaces, namely a
3-dimensional Einstein-Weyl geometry in Lorentzian signature. In
order to possess it an ODE must satisfy two point invariant
differential conditions on the function $F$, one of them being the
W\"unschmann condition. E. Cartan also showed that the Weyl
connection for this geometry may be immediately obtained from the
coframe $\theta^1,\ldots,\Omega_3$. Following Cartan's reasoning
S.-S. Chern studied third-order ODEs modulo contact
transformations and constructed a 10-dimensional bundle $\P\to
J^2$ equipped with a coframe
$\theta^1,\ldots\theta^4,\Omega_1,\ldots,\Omega_6$. Iff an ODE
satisfies the W\"unschmann condition then it has a
three-dimensional Lorentzian conformal geometry on the solution
space while the coframe becomes the $\o(3,2)$-valued normal
conformal connection for the geometry. In both these cases the
conformal metric is precisely the metric appearing implicitly in
K. W\"unschmann's thesis. Later H. Sato and A. Yoshikawa
\cite{Sat} applying N. Tanaka's theory \cite{Tan} constructed a
Cartan normal connection for arbitrary third-order ODEs (not only
of the W\"unschmann type) and showed how its curvature is
expressed by the contact relative invariants.

The Lorentzian geometry on the solution space was rediscovered
fifty years after S.-S. Chern from the perspective of General
Relativity. In a series of papers \cite{nsf2,nsf3,nsf4,New4} E.T.
Newman et al developed the Null Surface Formulation (NSF), an
alternate approach to General Relativity. In NSF the basic concept
is a family of hypersurfaces on a manifold $\M^4$; the
hypersurfaces are defined as level sets $Z(x^\mu,s,s^*)=const$ of
a real function\footnote{Here $s^*$ is the complex conjugate of
$s$.} $Z$ on $\M\times S^2$, where $(x^\mu)\in\M^4$, $s\in S^2$.
Starting from these data a Lorentzian {\em conformal} metric is
constructed by the property that these hypersurfaces are its null
hypersurfaces. The function $Z(x^\mu,s,s^*)$ is interpreted as the
general solution to a system of PDEs \ben\begin{aligned}
Z_{ss}=&\Lambda(s, s^*, Z, Z_s, Z_{s^*}, Z_{ss^*}), \notag
\\
Z_{s^*s^*}=&\Lambda^*(s, s^*, Z, Z_s, Z_{s^*}, Z_{ss^*}) \notag
\end{aligned}\een
for a real function $Z(s)$ of a complex variable $s$.
Consequently, coordinates $x^\mu$ are treated as constants of
integration, which turns $\M^4$ into the solution space of the
PDEs. One of main results of NSF was proving that a family of
hypersurfaces is a family of null hypersurfaces for a Lorentzian
conformal metric on $\M^4$ if and only if $\Lambda$ satisfies two
differential conditions, so called metricity conditions, which may
be viewed as a generalization of the W\"unschmann condition. It is
remarkable that three-dimensional version of the NSF leads
immediately to third-order ODEs, \cite{Tod,New2, New1, New3}. In
this case there is a one-parameter family of surfaces in $\M^3$
given by $Z(x^i,s)=const$, where $(x^i)\in\M^3$ and $s\in S^1$ are
real. $Z(x^i,s)$ is identified with the general solution of a
third order ODE $Z_{sss}=\Lambda(s,Z,Z_s,Z_{ss})$, $M^3$ is the
solution space and $\M^3\times S^1$ is identified with $\J^2$. The
construction of a conformal geometry from the family of null
surfaces is fully equivalent to Chern's construction and one
metricity condition obtained for $\Lambda$ in this case is
precisely the W\"unschmann condition.

Apart from the 3-dimensional Lorentzian conformal geometry, there
has been some interest in 3-dimensional Einstein-Weyl geometry,
mainly from the perspective of the theory of twistors and
integrable systems by N. Hitchin \cite{Hitch}, R. Ward
\cite{Ward}, C. LeBrun \cite{Leb} and P. Tod, M. Dunajski et al
\cite{Jt,Dun1,Dun2}; for discussion of the link between the
Einstein-Weyl spaces and third-order ODEs see \cite{Nur1} and
\cite{Tod}.

 P. Nurowski, following the ideas of E.
Cartan, proposed a programme of systematic study of geometries
related to differential equations, including second- and
third-order ODEs. In this programme,
\cite{New3,God3,Godode5,Nur1,Nur3}, both new and already known
geometries associated with differential equations are supposed to
be constructed by the Cartan equivalence method and are to be
characterized in the language of Cartan connections associated
with them. In particular in \cite{Nur1} new examples of geometries
associated with ordinary differential equations were given,
including a conformal geometry with special holonomy $G_2$ from
ODEs of the Monge type. Partial results on geometries of
third-order ODEs were given in \cite{Nur1, New3, God3} but the
full analysis of these geometries has not been published so far
and this paper aims to fill this gap.

Geometry of third-order ODEs is a part of broader issue of
geometry of differential equations in general. Regarding ODEs of
order two, we owe classical results including construction of
point invariants to S. Lie \cite{Lie} and M. Tresse \cite{Tre}. In
particular E. Cartan \cite{Car5} constructed a two-dimensional
projective differential geometry on the solution spaces of some
second-order ODEs. This geometry was further studied in \cite{NN}
and \cite{Nur3}, the latter paper pursues the analogy between
geometry of three-dimensional CR structures and second-order ODEs
and provides a construction of counterparts of the Fefferman
metrics for the ODEs. Classification of second-order ODEs
possessing Lie groups of fibre-preserving symmetries was done by
L. Hsu and N. Kamran \cite{Hsu}. Geometry on the solution space of
certain four-order ODEs (satisfying two differential conditions),
which is given by the four-dimensional irreducible representation
of $GL(2,\real)$ and has exotic $GL(2,\real)$ holonomy was
discovered and studied by R. Bryant \cite{Bry}, see also
\cite{Nur4}. The $GL(2,\real)$ geometry of fifth-order ODEs has
been recently studied by M. Godli\'nski and P. Nurowski
\cite{Godode5}.

The more general problems yet are existence and properties of
geometry on solution spaces of arbitrary ODEs. The problem of
existence was solved by B. Doubrov \cite{Dubgl}, who  proved that
an $n$th-order ODE $n\geq 3$, modulo contact transformations, has
 a geometry based on the irreducible
$n$-dimensional representation of $GL(2,\real)$ provided that it
satisfies  $n-2$ scalar differential conditions. An implicit
method of constructing these conditions was given in
\cite{Dubwil}. Properties of the $GL(2,\real)$ geometries of ODEs
are still an open problem; they were studied in \cite{Godode5},
where the Doubrov conditions were interpreted as higher order
counterparts of the W\"unschmann condition, and by M. Dunajski and
P. Tod \cite{Dun3}. %from the point of view paraconformal
%geometry \cite{Est}.

Almost all the above papers deal with geometries on solution
spaces but one can also consider other geometries, including those
defined on various jet spaces. The most general result on such
geometries \cite{Dub2} comes from T. Morimoto's nilpotent geometry
\cite{Mor1,Mor2}. It concludes that with a system of ODEs there is
associated a filtration on a suitable jet space together with a
canonical Cartan connection.

\subsection{Results of the paper}
\noindent We start from the equivalence problems formulated in
terms of $G$-structures \eqref{e.G_c} -- \eqref{e.G_f} and apply
Cartan's method to obtain manifolds $\P$ equipped with the
coframes encoding all the invariant information about ODEs. Next
we show how to read the principal bundle structures of these
manifolds over distinct bases. Usually it is the structure over
the solution space $\S$ which is the most interesting, but we also
consider structures over $\J^1$, $\J^2$ and certain
six-dimensional manifold $\M^6$, which appears naturally. When the
structure of a bundle is established then the invariant coframe
defines a Cartan connection on $\P\to base$, usually under
additional conditions playing similar role to the W\"unschmann
condition. In order to obtain the geometries and the
W\"unschmann-like conditions we often apply the method of
construction by Lie transport and projection. The most symmetric
cases when the connections are flat or have constant curvature
provide homogeneous models for the geometries. Since the dimension
of $\S$ and $\J^1$ is three and we build geometries with at least
two-dimensional structural group then the homogeneous models are
given by the ODEs with at least five-dimensional symmetry group.

\subsubsection*{Contact geometries} Section \ref{ch.contact} is
devoted to the geometries of the ODEs modulo contact
transformations. The only equations possessing at least five-
dimensional contact symmetry group are the linear equations with
constant coefficients, that is \ben y'''=0 \een with the symmetry
group $O(3,2)$ and \ben y'''=-2\mu y'+y,\qquad \mu\in \real, \een
mutually non-equivalent for distinct $\mu$, with the symmetry
group $\real^2\ltimes_{\mu}\real^3$. Sections \ref{s.c.th} to
\ref{s.c6d} discuss geometries whose homogeneous model is
generated by $y'''=0$. In section \ref{s.c.th} we state the main
theorem in that section, theorem \ref{th.c.1}, which describes the
geometry on $\J^2$. It may be recapitulated as follows.

\begin{theo}[Theorem \ref{th.c.1}] The contact invariant information
about an equation $y'''=F(x,y,y',y'')$ is given by the following
data
\begin{itemize}
\item[i)] The principal fibre bundle $H_6\to\P\to\J^2$, where
$\dim\P=10$ and $H_6$ is a six-dimensional subgroup of $O(3,2)$
\item[ii)] The coframe
$(\theta^1,\theta^2,\theta^3,\theta^4,\Omega_1,\Omega_2,\Omega_3,\Omega_4,\Omega_5,\Omega_6)$
on $\P$ which defines  the $\o(3,2)\cong\sp(4,\real)$ Cartan
normal connection $\wh{\omega}^c$ on $\P$.
\end{itemize}
The coframe and the connection $\wh{\omega}_c$ are given
explicitly in terms of $F$ and its derivatives. There are two
basic relative invariants for this geometry: the W\"unschmann
invariant and $F_{y''y''y''y''}$.\end{theo} This theorem is almost
identical to the result proved in \cite{Sat} and the only new
element we add here is the explicit formula for the connection.
For completeness, section \ref{s.c.proof} contains a proof of the
theorem, which follows S.-S. Chern's construction of the coframe
and the construction of the normal connection of \cite{Sat}.

Sections \ref{s.c.conf} to \ref{s.c6d} discuss next three
geometries generated by $\wh{\omega}^c$, those on $\J^1$, $\S$ and
certain six-dimensional manifold $\M^6$ which appears in a natural
way. These results are summarized as follows.
\begin{theo}
The connection $\wh{\omega}^c$ of theorem \ref{th.c.1} has
fourfold interpretation.
\begin{itemize}
 \item[1.] It is always the normal $\o(3,2)$ Cartan connection  on
 $\J^2$ (Chern-Sato-Yoshikawa construction.)
\item[2.] If the ODE has vanishing W\"unschmann condition then
$\wh{\omega}^c$ is the normal Lorentzian conformal connection  for
the Lorentzian structure on the solution space (Chern-NSF
construction.) \item[3.] If the ODE satisfies $F_{y''y''y''y''}=0$
then $\wh{\omega}^c$ becomes the normal Cartan connection for the
contact projective structure on $\J^1$ generated by the family of
solutions of the ODE. \item[4.] $\wh{\omega}^c$ is the
$\sp(4,\real)$-part of the $\o(4,4)$ normal conformal connection
for a six-dimensional split conformal geometry on $\M^6$ with
special holonomy $\sp(4,\real)\semi{.}\real^5$.
\end{itemize}
\end{theo}

In section \ref{s.c.furred} we turn to geometries, whose
homogeneous models are provided by the equation $y'''=-2\mu y'+y$.
Following Chern we reduce the bundle $\P$ to its five-dimensional
subbundle. Then we find that
\begin{theo}[Theorem \ref{cor.c.geom5d}]
Every ODE satisfying some contact invariant condition
$\inc{a}{}[F]=\mu=const$ has a $\real^2$ geometry on its solution
space together with a $\real^2$ linear connection from the
invariant coframe. The action of the algebra $\real^2$ on $\S$ is
given by \ben
   \bma
      u & v & 0 \\
      -\mu v & u & v \\
      v & -\mu v & u
   \ema.
 \een
\end{theo}
This geometry seems to be a generalization of Chern's `cone
geometry' which was associated with the equation $y'''=-y$ and
briefly mentioned to exist for arbitrary ODEs. In our construction
the action of $\real^2$ depends on the characteristic polynomial
of respective linear equation, and we get a real cone geometry
provided that it has three distinct roots.

\subsubsection*{Point geometries} In section \ref{ch.point}
we study the geometries associated with the ODEs modulo point
transformations. Sections \ref{s.p.th} to \ref{s.w6d} deal with
geometries modelled on $y'''=0$. Our approach is analogous to the
contact case and results are similar. They may be summarized as
follows.
\begin{theo}
The following statements hold
\begin{itemize}
\item[1.] The point invariant information about
$y'''=F(x,y,y',y'')$ is given by the seven-dimensional principal
bundle $H_3\to\P\to\J^2$ together with the coframe
$\hp{1},\hp{2},\hp{3},\hp{4},\vp{1},\vp{2},\vp{3}$ on $\P$, which
defines the $\co(2,1)\semi{.}\real^3$ Cartan connection
$\wh{\omega}^p$ (Cartan construction.) \item[2.] If the ODE has
vanishing W\"unschmann  and Cartan invariants then it has the
Einstein-Weyl geometry on $\S$ and the Weyl connection is given by
$\wh{\omega}^p$ (Cartan construction.) \item[3.] If the ODE
satisfies $F_{y''y''y''}=0$ then it has the point-projective
structure on $\J^1$ generated by the family of solutions of the
ODE. \item[4.] For any ODE there exists the split signature
six-dimensional Weyl geometry on certain manifold $\M^6$, which is
never Einstein.
\end{itemize}
\end{theo}

A new construction, which does not have a contact counterpart, is
considered in section \ref{s.lor3}. This is a Lorentzian {\em
metric structure} on the solution space $\S$. Its construction
follows immediately from the Einstein-Weyl geometry. If the Ricci
scalar of the Weyl connection is non-zero, then it is a weighted
conformal function and may be fixed to a constant by an
appropriate choice of the conformal gauge. The homogeneous models
of this geometry are associated with \ben
y'''=\frac{3\,y''^2}{2\,y'} \een if the Ricci scalar is negative,
and \ben y'''=\frac{3y''^2y'}{y'^2+1}\een if the Ricci scalar is
positive. Their point symmetry groups are $O(2,2)$ and $O(4)$
respectively. Both these equations are contact equivalent to
$y'''=0$.

\subsubsection*{Fibre-preserving geometries} Section \ref{ch.fp}
is devoted to the geometries of ODEs modulo fibre-preserving
transformations. We obtain a seven-dimensional bundle and the
$\co(2,1)\semi{.}\real^3$-valued Cartan connection $\wh{\omega}^f$
on it. Since both $\wh{\omega}^f$ and $\wh{\omega}^p$ of the point
case take value in the same algebra these cases are very similar
to each other. Indeed, we show that one can recover
$\wh{\omega}^f$ from $\wh{\omega}^p$ just by appending one
function on the bundle. As a consequence, the geometries of the
fibre-preserving case are obtained from their point counterparts
by appending the object generated by the function.

We did not study obvious or not interesting geometries. In the
point and fibre-preserving case geometries of $y'''=-2\mu y'+y$
are the same to what we have in the contact case, since the
respective symmetry groups are the same. Also the fibre-preserving
geometry on $\J^1$ does not seem to be worth studying.

To summarize, the following material contained in this work is
new: sections \ref{s.c-p} to \ref{s.c.furred} excluding theorem
\ref{th.c.2} and sections \ref{s.p-p} to \ref{s.fp}. Other
sections contain a reformulation and an extension of already known
results.

All our calculations were performed or checked using the symbolic
calculations program Maple.

\subsection{Notation}\label{s.notation} In what follows we use the following symbols,
in particular $W$ denotes the W\"unschmann invariant.
\begin{align}
 F=&F(x,y,p,q), \notag \\
\D=&\partial_x+p\partial_y+q\partial_p+F\partial_q, \\ % \label{e.defD} \\
 K=&\tfrac{1}{6}\D F_q -\tfrac{1}{9}F_q^2-\tfrac{1}{2}F_p, \label{e.defK} \\
 L=&\tfrac13F_{qq}K-\tfrac13F_qK_q-K_p-\tfrac13F_{qy},  \\ % \label{e.defL} \\
 M=&2K_{qq}K-2K_{qy}+\tfrac13F_{qq}L-\tfrac23F_qL_q-2L_p,  \\ % \label{e.defM} \\
 W=&\left(\D-\tfrac{2}{3}F_q\right)K+F_y, \label{e.defW} \\
 Z=&\frac{\D W}{W} -F_q.  % \label{e.defZ}
\end{align}
Parentheses denote sets of objects: $(a_1,\ldots,a_k)$ is the set
consisting of $a_1,\ldots,a_k$. In particular this symbol denotes
bases of vector spaces as well as coordinate systems, frames and
coframes on manifolds.  The linear span of vectors or covectors
$a_1,\ldots,a_k$ is denoted by $<a_1,\ldots,a_k>$. If
$a_1,\ldots,a_k$ are vector fields or one-forms on a manifold,
then the above symbol denotes the distribution or the simple ideal
generated by them. The symmetric tensor product of two one-forms
or vector fields $\alpha$ and $\beta$ is denoted by
$\alpha\beta=\tfrac12(\alpha\otimes\beta+\beta\otimes\alpha)$. The
symbols $A_{(\mu\nu)}$ and $A_{[\mu\nu]}$ denote symmetrization
and antisymmetrization of a tensor $A_{\mu\nu}$ respectively. For
a metric $g$ of signature $(k,l)$ the group $CO(k,l)$ is defined
to be
$$CO(k,l)=\{A\in GL(k+l,\real)\,|\quad A^TgA=e^\lambda g,\quad\lambda\in\real \}. $$
Its Lie algebra
$$\co(k,l)=\{a\in \gl(k+l,\real)\,|\quad a^Tg+ga=\lambda g,\quad\lambda\in\real \}. $$
A semidirect product of two Lie groups $G$ and $H$, where $G$ acts
on $H$ is denoted by $ G\ltimes H$. A semidirect product of their
Lie algebras is denoted by $\g\semi{.}\h$. If the action depends
of a parameter $\mu$ then we add an subscript: $\ltimes_\mu$ and
$\semi{\mu}$.

\section{Towards Cartan connections}
\label{ch.prel}

\subsection{ODEs as G-structures on $\J^2$} Following Cartan and Chern
we begin with the space $\J^2$ of second jets of curves in
$\real^2$ (see \cite{Olv1} for description of jet spaces) with
coordinate system $(x,y,p,q)$ where $p$ and $q$ denote first and
second derivative $y'$ and $y''$ for a curve $x\mapsto(x,y(x))$ in
$\real^2$, so that this curve lifts to a curve
$x\mapsto(x,y(x),y'(x),y''(x))$ in $\J^2$. Any solution $y=f(x)$
of $y'''=F(x,y,y',y'')$ is uniquely defined by a choice of
$f(x_0)$, $f'(x_0)$ and $f''(x_0)$ at some $x_0$. Since that
choice is equivalent to a choice of a point in $\J^2$, there
passes exactly one solution $(x,f(x),f'(x),f''(x))$ through any
point of $\J^2$. Therefore the solutions form a (local) congruence
on $\J^2$, which can be described by its annihilating simple
ideal. Let us choose a coframe $(\omega^i)$ on $\J^2$:
\be\label{e.omega}\begin{aligned}
 \omega^1&=\der y -p\der x,  \\
 \omega^2&=\der p- q\der x,  \\
 \omega^3&=\der q- F(x,y,p,q)\der x,  \\
 \omega^4&=\der x.
\end{aligned}
\ee Each solution $y=f(x)$  is fully described by the two
conditions: forms $\omega^1$, $\omega^2$, $\omega^3$ vanish on the
curve $t\mapsto(t,f(t),f'(t),f''(t))$ and, since this defines a
solution modulo transformations of $x$, $\omega^4=\der t$ on this
curve.

Suppose now that a equation $y'''=F(x,y,,y',y'')$ undergoes a
contact, point or fibre-preserving transformation. Then
\eqref{e.omega} transform by \be\bal
 \omega^1&\mapsto \bar{\omega}^1=u_1\omega^1, \\
 \omega^2&\mapsto \bar{\omega}^2=u_2\omega^1+u_3\omega^2,  \\
 \omega^3&\mapsto \bar{\omega}^3=u_4\omega^1+u_5\omega^2+u_6\omega^3,\\
 \omega^4&\mapsto \bar{\omega}^4=u_8\omega^1+u_9\omega^2+u_7\omega^4, \\
\eal\label{e.trans_kor} \ee with some functions $u_1,\ldots,u_9$
defined on $\J^2$ and determined by a particular choice of
transformation, for instance \begin{eqnarray*}
&u_1=\phi_y-\psi\chi_y,\\
&u_7=\Der \chi,\quad u_8=\chi_y,  \quad u_9= \chi_p.
\end{eqnarray*}
In particular, $u_9=0$ in the point case and $u_8=u_9=0$ in the
fibre-preserving case. Since the transformations are
non-degenerate, the condition $u_1u_3u_5u_7\neq 0$ is always
satisfied and transformations \eqref{e.trans_kor} form local
pseudogroups. Thus we have. \begin{lemma} A class of contact
equivalent third-order ODEs is a local G-structure\footnote{Here
we use local trivializations of G-structures.} $G_c\times\J^2$
defined by the property that the coframe
$(\omega^1,\omega^2,\omega^3,\omega^4)$ belongs to it and the
structure group is given by \be\label{e.G_c}
G_c= \bma u_1 & 0 & 0 & 0 \\ u_2 & u_3 & 0 & 0 \\
u_4 & u_5 & u_6& 0\\
u_8 & u_9 & 0 & u_7 \ema. \ee For a class of point equivalent ODEs
the structure group $G_p\subset G_c$ is given by \be\label{e.G_p}
G_p= \bma u_1 & 0 & 0 & 0 \\ u_2 & u_3 & 0 & 0 \\
u_4 & u_5 & u_6& 0\\
u_8 & 0 & 0 & u_7 \ema, \ee whereas for a class of
fibre-preserving equivalent ODEs the structure group $G_f\subset
G_p$ is given by \be\label{e.G_f}
G_f= \bma u_1 & 0 & 0 & 0 \\ u_2 & u_3 & 0 & 0 \\
u_4 & u_5 & u_6& 0\\
0 & 0 & 0 & u_7 \ema. \ee
\end{lemma}

\subsection{Cartan connections}
The version of Cartan's equivalence method \cite{Car4} we employ
below is explained in books by R. Gardner \cite{Gar} and P. Olver
\cite{Olv2}. Some its aspects are also discussed by S. Sternberg
\cite{Ste} and S. Kobayashi \cite{Kob1}. Its idea is the
following: starting from $G_c\times\J^2$ (or $G_p\times\J^2$ or
$G_f\times\J^2$) one constructs a new principal bundle
$H\to\P\to\J^2$ (with a new group $H$) equipped with one fixed
coframe $(\theta^i,\Omega_{\mu})$ built in some geometric way and
such that it encodes all the local invariant information about
$G_c\times\J^2$ (or $G_p\times\J^2$ or $G_f\times\J^2$
respectively) through its structural equations. This leads to the
notion of Cartan  connection, defined here after \cite{Kob1}.

\begin{definition}
Let $H\to\P\to\M$ be a principal bundle and let $G$ be a Lie group
such that $H$ is its closed subgroup and $\dim G=\dim\P$. A Cartan
connection of type $(G,H)$ on $\P$ is a one-form $\wh{\omega}$
taking values in the Lie algebra $\goth{g}$ of $G$ and satisfying
the following conditions:
\begin{itemize}
\item[i)] $\wh{\omega}_u:T_u\P\to\goth{g}$ for every $u\in\P$ is
an isomorphism of vector spaces \item[ii)]
$A^*\hook\,\wh{\omega}=A$ for every $A\in\goth{h}$ and the
corresponding fundamental field $A^*$ \item[iii)]
$R^*_h\wh{\omega}=\Ad(h^{-1})\wh{\omega}$ for $h\in H$.
\end{itemize}
\end{definition}

The curvature of a Cartan connection is defined as follows \ben
\wh{K}(X,Y)=\der\wh{\omega}(X,Y)+\frac{1}{2}[\wh{\omega}(X),\wh{\omega}(Y)].
\een The curvature is horizontal which means that it vanishes on
each fundamental vector field $A^*$: \ben A^*\hook\, \wh{K}=0.\een
Horizontality of the curvature is locally equivalent with property
iii) of the connection. Cartan connections with vanishing
curvature are called flat.

\begin{example}
A standard example of Cartan connection is obtained by taking the
principal fibre bundle of a Lie group $\P=G$ over its homogeneous
space $\M=G/H$ and by defining $\wh{\omega}$ to be the
Maurer-Cartan form $\wh{\omega}=g^{-1}\der g$, $g\in G$. So
defined $\wh{\omega}$ is a flat Cartan connection on $\P$.
\end{example}

In general, the curvature does not have to vanish and a Cartan
connection is an object that generalizes the notion of the
Maurer-Cartan form on a Lie group.

Two 3rd order ODEs $y'''=F(x,y,y',y'')$ and
$y'''=\bar{F}(x,y,y',y'')$ are contact/point/fibre-preserving
equivalent if and only if their associated Cartan connections are
diffeomorphic, that is there exists a local bundle diffeomorphism
$\Phi\colon\bar{\P}\supset\bar{\O}\to\O\subset\P$ such that
$\Phi^*\wh{\omega}=\overline{\wh{\omega}}$. Technically speaking
curvature coefficients are rational functions of vertical bundle
coordinates $(u_\mu)$ and relative invariants defined as follows.
\begin{definition}
A contact/point/fibre-preserving relative invariant $I[F]$ of
$y'''=F(x,y,y',y'')$ is a function of $F$ and its derivatives such
that if $I[F]\equiv 0$ for an equation $y'''=F(x,y,y',y'')$ then
$I[\bar{F}]\equiv 0$ for every equation $y'''=\bar{F}(x,y,y',y'')$
contact/point/fibre-preserving equivalent to it.
\end{definition}
From relative invariants contained in the curvature of the
associated Cartan connection one may recover the full set of
invariants by taking consecutive coframe derivatives.

In theorem \ref{th.c.1} and proposition \ref{prop.c-p} we need a
concept of normal Cartan connection in the sense of N. Tanaka. We
briefly remind this notion below. For details of Tanaka's theory
see \cite{Tan2,Tan}.
\begin{definition}
A semisimple Lie algebra $\g$ is graded if it has a vector space
decomposition \ben
\g=\g_{-k}\oplus\ldots\oplus\g_{-1}\oplus\g_0\oplus\g_1\oplus\ldots\oplus\g_{k}
\een such that \ben[\g_i,\g_j]\subset \g_{i+j}\een and
$\g_{-k}\oplus\ldots\oplus\g_{-1}$ is generated by $\g_{-1}$.
\end{definition}

Let us suppose that $\g$ is a semisimple graded Lie algebra and
denote $\m=\g_{-k}\oplus\ldots\oplus\g_{-1}$,
$\h=\g_0\oplus\ldots\oplus\g_k$. Let us consider a $\g$ valued
Cartan connection $\wh{\omega}$ on a bundle $H\to\P\to\M$, where
the Lie algebra of $H$ is $\h$. Fix a point $p\in\P$. The
decomposition $\g=\m\oplus\h$ defines in $T_p\P$ the complement
$\H_p$ of the vertical space $\V_p$. Therefore we have
$T_p\P=\V_p\oplus\H_p$, $\wh{\omega}(\V_p)=\h$ and
$\wh{\omega}(\H_p)=\m$. The curvature $\wh{K}_p=$
$(\der\wh{\omega}+\wh{\omega}\w\wh{\omega})_p$ at $p$ is then
characterized by the tensor $\kappa_p\in\Hom(\wedge^2\m,\g)$ given
by
 \be\label{e.c.kappa}\kappa_p(A,B)=\wh{K}_p(\wh{\omega}^{-1}_p(A), \wh{\omega}^{-1}_p(B)),\quad
 A,B\in\m.\ee
The function $\kappa\colon\P\to\Hom(\wedge^2\m,\g)$ is called the
structure function.

In the space $\Hom(\wedge^2\m,\g)$ let us define
$\Hom^1(\wedge^2\m,\g)$ to be the space of all
$\alpha\in\Hom(\wedge^2\m,\g)$ fulfilling
\ben\alpha(\g_i,\g_j)\subset \g_{i+j+1}\oplus\ldots\oplus\g_k
\quad\text{for}\quad i,j<0.\een Since the Killing form $B$ of $\g$
is non-degenerate and satisfies $B(\g_p,\g_q)=0$ for $p\neq -q$,
one can identify $\m^*$ with $\g_1\oplus\ldots\oplus\g_k$. For a
basis $(e_1,\ldots,e_m)$ of $\m$ let $(e^*_1,\ldots,e^*_m)$ denote
the unique basis of $\g_1\oplus\ldots\oplus\g_k$ such that
$B(e_i,e^*_j)=\delta_{ij}$. Tanaka considered the following
complex \ben
\ldots\longrightarrow\Hom(\wedge^{q+1}\m,\g)\overset{\partial^*}{\longrightarrow}\Hom(\wedge^q\m,\g)\longrightarrow\ldots
\een with
$\partial^*\colon\Hom(\wedge^{q+1}\m,\g)\to\Hom(\wedge^q\m,\g)$
given by the following formula \ben
\begin{aligned}
(\partial^*\alpha)&(A_1\w\ldots\w A_q)=\sum_i[e^*_i,\alpha(e_i\w
A_1\w\ldots\w A_q)]
 \\&+\tfrac{1}{2}\sum_{i,j}\alpha([e^*_j,A_i]_\m\w e_j\w A_1\w\ldots\w\hat{A_i}\w\ldots\w A_q),
\end{aligned}
\een where $\alpha\in\Hom(\wedge^q\m,\g)$, $A_1,\ldots\,A_q\in\m$,
$(e_i)$ is any basis in $\m$ and $[\,\, ,\,]_\m$ denotes the
$\m$-component of the bracket with respect to the decomposition
$\g=\m\oplus\h$. Finally, N. Tanaka \cite{Tan} introduced the
notion of normal connection, the definition below is given in the
language of \cite{Cap}.
\begin{definition}\label{def.Tanakanorm}
A Cartan connection $\wh{\omega}$ as above is normal if its
structure function $\kappa$ fulfills the following conditions \ben
\begin{aligned}
\text{i)}& &\qquad &\kappa\in\Hom^1(\wedge^2\m,\g),\\
\text{ii)}& &\qquad &\partial^*\kappa=0.
\end{aligned}
\een
\end{definition}

\section{Geometries of ODEs modulo contact transformations of
variables}\label{ch.contact}

\subsection{Cartan connection on ten-dimensional
bundle}\label{s.c.th} \noindent We formulate a theorem about the
main structure which is associated with third-order ODEs modulo
contact transformations of variables, an $\sp(4,\real)$ Cartan
connection on the bundle $\P^c\to\J^2$. This structure will serve
as a starting point for further analyzing of geometries of ODEs.

\begin{theorem}\label{th.c.1}
To every third order ODE $y'''=F(x,y,y',y'')$ there are associated
the following data.
\begin{itemize}
\item[i)] The principal fibre bundle $H_6\to\P^c\to\J^2$, where
$\dim\P^c=10$ and $H_6$ is the following six-dimensional subgroup
of $SP(4,\real)$ \be \label{e.c.H6}  H_6=\bma \sqrt{u_1}, &
\frac12\frac{u_2}{\sqrt{u_1}}, & -\frac12\frac{u_4}{\sqrt{u_1}}, &
 \tfrac{1}{24}\tfrac{u_2^2u_5}{u_1^{3/2}u_3}-\tfrac12\sqrt{u_1}\,u_6 \\\\
 0 & \tfrac{u_3}{\sqrt{u_1}}, & -\tfrac{u_5}{\sqrt{u_1}},
 & \tfrac12\tfrac{u_2u_5-u_3u_4}{u_1^{3/2}} \\\\
 0 & 0 & \tfrac{\sqrt{u_1}}{u_3}, &
-\tfrac12\tfrac{u_2}{\sqrt{u_1}\,\u_3} \\\\
  0 & 0 & 0 & \tfrac{1}{\sqrt{u_1}}\ema.
\ee \item[ii)] The coframe
$(\theta^1,\theta^2,\theta^3,\theta^4,\Omega_1,\Omega_2,\Omega_3,\Omega_4,\Omega_5,\Omega_6)$
on $\P^c$, which defines the $\sp(4,\real)$-valued Cartan normal
connection $\wh{\omega}^c$ on $\P^c$ by \be\label{e.c.conn_sp}
 \wh{\omega}^c=\bma \tfrac{1}{2}\vc{1} & \tfrac{1}{2}\vc{2} & -\tfrac{1}{2}\vc{4} & -\tfrac{1}{4}\vc{6} \\\\
             \hc{4} & \vc{3}-\tfrac{1}{2}\vc{1} & -\vc{5} & -\tfrac{1}{2}\vc{4} \\\\
             \hc{2} & \hc{3} & \tfrac{1}{2}\vc{1}-\vc{3} & -\tfrac{1}{2}\vc{2} \\\\
             2\hc{1} & \hc{2} & -\hc{4} & -\tfrac{1}{2}\vc{1}
        \ema.
\ee
\end{itemize}

Two 3rd order ODEs $y'''=F(x,y,y',y'')$ and
$y'''=\bar{F}(x,y,y',y'')$ are locally contact equivalent if and
only if their associated Cartan connections are locally
diffeomorphic, that is there exists a local bundle diffeomorphism
$\Phi\colon\bar{\P}^c\supset\bar{\O}\to\O\subset\P^c$ such that
$$\Phi^*\wh{\omega}^c=\overline{\wh{\omega}^c}. $$

The connection $\wh{\omega}^c$ has the following explicit form.
Let $(x,$ $y,$ $p,$ $q,$ $u_1,$ $u_2,$ $u_3,$ $u_4,$ $u_5,$
$u_6)$, $(x^i,u_\mu)$ for short, be a local coordinate system in
$\P^c$, which is compatible with the local trivialization
$\P^c=H_6\times\J^2$, that is $(x^i)=(x,y,p,q)$ are coordinates in
$\J^2$ and $(u_\mu)$ are coordinates in $H_6$ as in
\eqref{e.c.H6}. Then the value of $\wh{\omega}^c$ at the point
$(x^i,u_\mu)$ in $\P^c$ is given by \ben
\wh{\omega}^c(x^i,u_\mu)=u^{-1}\,\omega^c\,u+u^{-1}\der u \een
where $u$ denotes the matrix \eqref{e.c.H6} and \ben {\omega}^c=
\bma
\tfrac{1}{2}\vc{1}^0&\tfrac{1}{2}\vc{2}^0 & -\tfrac{1}{2}\vc{4}^0 & -\tfrac{1}{4}\vc{6}^0 \\\\
\omega^4 & \vc{3}^0-\tfrac{1}{2}\vc{1}^0 & -\vc{5}^0 & -\tfrac{1}{2}\vc{4}^0 \\\\
\omega^2 & \wt{\omega}^3 & \tfrac{1}{2}\vc{1}^0-\vc{3}^0 & -\tfrac{1}{2}\vc{2}^0 \\\\
2\omega^1 & \omega^2 & -\omega^4 & -\tfrac{1}{2}\vc{1}^0 \ema\een
is the connection $\wh{\omega}^c$ calculated at the point
$(x^i,u_1=1,u_2=0,u_3=1,u_4=0,u_5=0,u_6=0)$.
 The forms $\omega^1,\omega^2,\wt{\omega}^3,\omega^4$ read
\be\label{e.c.om}\begin{aligned}
\omega^1=&\der y-p\der x, \\
\omega^2=&\der p-q\der x, \\
\wt{\omega}^3=&\der q-F\der x-\tfrac13F_q(\der p-q\der x)+K(\der y-p\der x), \\
\omega^4=&\der x.
\end{aligned}
\ee The forms $\vc{1}^0,\ldots,\vc{6}^0$ read
\be\label{e.c.Om0}\begin{aligned}
\vc{1}^0=&-K_q\,\omega^1,  \\
\vc{2}^0=&\left(\tfrac13 W_q+L\right)\,\omega^1-K_q\,\omega^2-K\omega^4,  \\
\vc{3}^0=&-K_q\,\omega^1+\tfrac16F_{qq}\,\omega^2+\tfrac13F_q\omega^4, \\
\vc{4}^0=&-(\tfrac13W_{qq}+L_q)\,\omega^1+\tfrac12K_{qq}\,\omega^2,
\\
\vc{5}^0=&\tfrac12K_{qq}\,\omega^1-\tfrac16F_{qqq}\,\omega^2-\tfrac16F_{qq}\,\omega^4,
 \\
\vc{6}^0=&(\tfrac13\D(W_{qq})-\tfrac43W_{qp} -\tfrac13F_qW_{qq}
+\tfrac13F_{qqq}W+M)\,\omega^1+ \\
&+\tfrac13(F_{qqy}-F_{qqq}K-W_{qq})\,\omega^2-K_{qq}\,\wt{\omega}^3+ \\
&+(\tfrac23F_{qy}-\tfrac13F_{qq}K-2L-\tfrac43W_q)\,\omega^4.
\end{aligned}
\ee
\end{theorem}

\subsection{Proof of theorem 3.1} \label{s.c.proof} \noindent We
prove theorem \ref{th.c.1} by repeating Chern's construction
supplemented later in \cite{Sat}.

On the bundle $G_c\times\J^2$ there are four fixed, well defined
one-forms $(\theta^1,\theta^2,\theta^3,\theta^4)$, the components
of the canonical $\real^4$-valued form $\theta$ existing on the
frame bundle of $\J^2$. Let $(x)$ denote $x,y,p,q$ and $(g)$ be
coordinates in $G_c$ given by \eqref{e.G_c}. Let us choose a
coordinate system $(x,g)$ on $G_c\times\J^2$ compatible with the
local trivialization. Then $\theta^i$ at the point $(x,g^{-1})$
read \be \label{e.i.canon}\bal
 \theta^1=&u_1\omega^1, \\
 \theta^2=&u_2\omega^1+u_3\omega^2,  \\
 \theta^3=&u_4\omega^1+u_5\omega^2+u_6\omega^3,\\
 \theta^4=&u_8\omega^1+u_9\omega^2+u_7\omega^4. \\
\eal \ee

We seek a bundle on which $\theta^i$ are supplemented to a coframe
by certain new one-forms $\Omega_\mu$ chosen in a well-defined
geometric manner.

\paragraph*{Step 1} We calculate the exterior derivatives of $\theta^i$
on $G_c\times\J^2$ \be\label{e.c.red10}\begin{aligned}
 \der\theta^1=&\alpha_1\w\theta^1+T^1_{~jk}\theta^j\w\theta^k, \\
 \der\theta^2=&\alpha_2\w\theta^1+\alpha_3\w\theta^2+T^2_{~jk}\theta^j\w\theta^k, \\
 \der\theta^3=&\alpha_4\w\theta^1+\alpha_5\w\theta^2+\alpha_6\w\theta^3+T^3_{~jk}\theta^j\w\theta^k, \\
 \der\theta^4=&\alpha_8\w\theta^1+\alpha_9\w\theta^2+ \alpha_7\w\theta^4
 +T^4_{~jk}\theta^j\w\theta^k,
 \end{aligned} \ee
where $\alpha_\mu$ are the entries of the matrix $\der
g^i_{~k}\cdot g^{-1k}_{~j}$ and $T^i_{~jk}$ are some functions on
$G_c\times\J^2$. Next we collect $T^i_{~jk}\theta^j\w\theta^k$
terms
\begin{align}
 \der\theta^1=&\left(\alpha_1-T^1_{~12}\theta^2-T^1_{~13}\theta^3-T^1_{~14}\theta^4\right)\w\theta^1\nonumber \\
  &+T^1_{~23}\theta^2\w\theta^3+T^1_{~24}\theta^2\w\theta^4+T^1_{~34}\theta^3\w\theta^4,\nonumber \\
 \der\theta^2=&\left(\alpha_2-T^2_{~12}\theta^2-T^2_{~13}\theta^3-T^2_{~14}\theta^4\right)\w\theta^1 \nonumber\\
   &+\left(\alpha_3-T^2_{~23}\theta^3-T^2_{~24}\theta^4\right)\w\theta^2+T^2_{~34}\theta^3\w\theta^4,\label{e.c.red20} \\
 \der\theta^3=&\left(\alpha_4-T^3_{~12}\theta^2-T^3_{~13}\theta^3-T^3_{~14}\theta^4\right)\w\theta^1,\nonumber \\
   &+\left(\alpha_5-T^3_{~23}\theta^3-T^3_{~24}\theta^4\right)\w\theta^2
   +\left(\alpha_6-T^3_{~34}\theta^4\right)\w\theta^3\nonumber\\
 \der\theta^4=&\left(\alpha_8-T^4_{~12}\theta^2-T^4_{~13}\theta^3-T^4_{~14}\theta^4\right)\w\theta^1,\nonumber \\
   &+\left(\alpha_9-T^4_{~23}\theta^3-T^4_{~24}\theta^4\right)\w\theta^2
   +\left(\alpha_7+T^4_{~34}\theta^3\right)\w\theta^4\nonumber
 \end{align}
and introduce new 1-forms $\pi_\mu$ substituting the collected terms.
%$\pi_1=\alpha_1-T^1_{~12}\theta^2-T^1_{~13}\theta^3-T^1_{~14}\theta^4$, $\pi_2=\alpha_2-T^2_{~23}\theta^3-T^2_{~24}\theta^4$
%and so on.
Eq. \eqref{e.c.red10} now read \be\label{e.c.red30}\begin{aligned}
 \der\theta^1&=\pi_1\w\theta^1+\frac{u_1}{u_3 u_7}\theta^4\w\theta^2, \\
 \der\theta^2&=\pi_2\w\theta^1+\pi_3\w\theta^2+\frac{u_3}{u_6u_7}\theta^4\w\theta^3, \\
 \der\theta^3&=\pi_4\w\theta^1+\pi_5\w\theta^2+\pi_6\w\theta^3,\\
 \der\theta^4&=\pi_8\w\theta^1+\pi_9\w\theta^2+\pi_7\w\theta^4,
\end{aligned}\ee
since $T^1_{~23}=T^1_{~34}=0$, $T^1_{~24}=-u_3u_7/u_1$ and
$T^2_{~34}=-u_3/(u_6u_7)$.

 The equations \eqref{e.c.red30} resemble structural equations
for a linear connection very much, however, here $\pi$ is not a
linear connection since it does not transform according to
$R^*_u\pi=\Ad(u^{-1})\pi$ along fibres of $G_c\times\J^2$. We may
think of $\pi$ as a connection in a broader meaning, that is a
horizontal distribution on $G_c\times\J^2$, which is not
necessarily right-invariant. Keeping this in mind we will refer to
$T^i_{~jk}$ as torsion. Thus $\pi$ is a connection chosen by the
demand that its torsion is `minimal', i.e. possesses as few terms
as possible. The forms $\pi_\mu$, which are candidates for the
sought forms $\Omega_\mu$, are not uniquely defined by equations
\eqref{e.c.red30}, for example the gauge $\pi_1\to\pi_1+f\theta^1$
leaves \eqref{e.c.red30} unchanged. Therefore our connection is
not uniquely defined by its torsion.

\paragraph*{Step 2} We reduce the bundle $G_c\times\J^2$. We
choose its subbundle, say $\P^{(1)}$, characterized by the
property that the torsion coefficients are constant on it. We
choose $\P^{(1)}$ such that $T^1_{~24}=-1,\,T^2_{~34}=-1$ on it.
Thus $\P^{(1)}$ is defined by \be\label{e.c.red_u6u7}
  u_6=\frac{u_3^2}{u_1},\quad\quad\quad u_7=\frac{u_1}{u_3}.
\ee  It is known \cite{Gar,Ste} that such a reduction preserves
the equivalence, in other words, two bundles are equivalent if and
only if their respective reductions are. Here $\P^{(1)}$ has the
seven-dimensional structural group \ben G_c^{(1)}=\bma u_1 & 0 & 0
& 0
\\ u_2 & u_3 & 0 & 0 \\ u_4 & u_5 & \tfrac{u_3^2}{u_1} & 0 \\ u_8 & u_9 & 0
& \tfrac{u_1}{u_3} \ema. \een

\paragraph*{Step 3} Next we pull-back $\theta^i$ and $\pi_\mu$ to $\P^{(1)}$.
But the new structural group $G^{(1)}_c$ is a seven-dimensional
subgroup of $G_c$, so
$(\theta^1,\ldots,\theta^4,\pi_1,\ldots,\pi_9)$ of
\eqref{e.c.red30} is not a coframe on $\P^{(1)}$ any longer, since
\ben \pi_6=2\pi_3-\pi_1 \mod(\theta^i), \quad\quad\quad
\pi_7=\pi_1-\pi_3
 \mod(\theta^i). \een
Taking this into account we recalculate \eqref{e.c.red30} and
gather the torsion terms. We choose the new connection
\ben\pi=\bma \pi_1 & 0 & 0 & 0 \\ \pi_2 & \pi_3 & 0 & 0 \\
\pi_4 & \pi_5 & 2\pi_3-\pi_1& 0\\
\pi_8 & \pi_9 & 0 & \pi_1-\pi_3 \ema \een so that its torsion is
minimal again. \be\label{e.c.red40}\begin{aligned}
 \der\theta^1&=\pi_1\w\theta^1+\theta^4\w\theta^2, \\
 \der\theta^2&=\pi_2\w\theta^1+\pi_3\w\theta^2+\theta^4\w\theta^3,\\
 \der\theta^3&=\pi_4\w\theta^1+\pi_5\w\theta^2+(2\pi_3-\pi_1)\w\theta^3
 +\left(\frac{3u_5}{u_3}-\frac{3u_2-u_3F_q}{u_1}\right)\theta^4\w\theta^3,\\
 \der\theta^4&=\pi_8\w\theta^1+\pi_9\w\theta^2+(\pi_1-\pi_3)\w\theta^4.
 \end{aligned}\ee

\paragraph*{Step 4} We repeat the steps 2. and 3. Firstly we
reduce $\P^{(1)}$ to the subbundle $\P^{(2)}\subset\P^{(1)}$
defined by the property that the only non-constant torsion
coefficient $T^3_{~34}$ in \eqref{e.c.red40} vanishes on it,
 \be\label{e.c.red_u5}
 u_5=\frac{u_3}{u_1}\left(u_2-\frac{1}{3}u_3F_q\right).
\ee Next we recalculate connection, re-collect the torsion and
make another reduction through the constant torsion condition ($K$
 is defined in \eqref{e.defK}.) \be\label{e.c.red_u4}
 u_4=\frac{u^2_3}{u_1}K+\frac{u_2^2}{2u_1}.
\ee

At this stage we have reduced the frame bundle $G_c\times\J^2$ to
the nine-dimensional subbundle $\P^{(3)}\to\J^2$, such that its
structural group is the following
\ben G^{(3)}_c=\bma u_1 & 0 & 0 & 0 \\
u_2 & u_3 & 0 & 0 \\
\tfrac{u_2^2}{u_1} & \tfrac{u_2u_3}{u_1} & \tfrac{u_3^3}{u_1} & 0 \\
u_8 & u_9 & 0 & \frac{u_1}{u_3} \ema \een and the frame dual to
$(\omega^1,\omega^2,\omega^3-\tfrac13F_q\omega^2+K\omega^1,\omega^4)$
belongs to $\P^{(3)}$. The structural equations on $\P^{(3)}$ read
after collecting \be\label{e.c.red50}\begin{aligned}
 \der\theta^1&=\pi_1\w\theta^1+\theta^4\w\theta^2, \\
 \der\theta^2&=\pi_2\w\theta^1+\pi_3\w\theta^2+\theta^4\w\theta^3, \\
 \der\theta^3&=\pi_2\w\theta^2+\left(2\pi_3-\pi_1\right)\w\theta^3
 +\frac{u_3^3}{u_1^3}W\theta^4\w\theta^1,\\
 \der\theta^4&=\pi_8\w\theta^1+\pi_9\w\theta^2+(\pi_1-\pi_3)\w\theta^4
\end{aligned}\ee
with some one-forms $\pi_1,\pi_2,\pi_3,\pi_8,\pi_9$. The function
$W$, defined in \eqref{e.defW}, is the W\"unsch\-mann invariant.
Thereby, as Chern observed, third-order ODEs fall into two main
contact inequivalent branches: the ODEs satisfying $W\neq 0$, and
those satisfying $W=0$.

Equations \eqref{e.c.red50} do not still define the forms
$\pi_\mu$ uniquely but only modulo the following transformations
\begin{align}
 \pi_1&\to  \pi_1+ 2t_1\theta^1,\nonumber \\
 \pi_2&\to  \pi_2+t_1\theta^2,\nonumber \\
 \pi_3&\to  \pi_3+t_1\theta^1, \label{e.c.prol10}\\
 \pi_8&\to  \pi_8+t_2\theta^1+t_3\theta^2+t_1\theta^4,\nonumber \\
 \pi_9&\to  \pi_9+t_3\theta^1+t_4\theta^2.\nonumber
\end{align}

At this point, there is no pattern of further reduction. If $W=0$
there are only constant torsion coefficients in \eqref{e.c.red50}
and we do not have any conditions to define a subbundle of
$\P^{(3)}$. In these circumstances we prolong $\P^{(3)}$.

\paragraph*{Step 5} Prolongation. On $\P^{(3)}$
there is no fixed coframe but only the coframe
$(\theta^1,\theta^2,\theta^3,\theta^4,\pi_1,\pi_2,\pi_3,\pi_8,\pi_9)$
given modulo \eqref{e.c.prol10}. But `a coframe given modulo $G$'
is a $G$-structure on $\P^{(3)}$ and we can deal with this new
structure on $P^{(3)}$ by means of the Cartan method. Let us
consider the bundle $G^{prol}\times\P^{(3)}$ then, where  \ben
G^{prol}=\bma 1 & 0 \\ t & 1 \ema \een reflects the freedom
\eqref{e.c.prol10} so that the block $t$ reads
\ben \bma 2t_1 & 0 & 0 & 0 \\ 0 & t_1 & 0 & 0 \\ t_1 & 0 & 0 & 0 \\
t_2 & t_3 & 0 & t_1 \\ t_3 & t_4 & 0 & 0 \ema. \een On
$\P^{(3)}\times G^{prol}$ there exist nine fixed one-forms
$\theta^1,\theta^2,\theta^3,\theta^4,\Pi_1,\Pi_2,\Pi_3,\Pi_8,\Pi_9$,
given by
\ben \bma \theta^i \\ \Pi_\mu \ema = \bma 1 & 0 \\
t & 1 \ema \bma \theta^i \\ \pi_\mu \ema, \een which is the
canonical one-form on $G^{prol}\times\P^{(3)}\to \P^{(3)}$.

\paragraph*{Step 6} Now we apply the method of reductions to the above
structure on $G^{prol}\times \P^{(3)}$. We calculate the exterior
derivatives of $(\theta^i,\Pi_\mu)$. The derivatives of $\theta^i$
take the form of \eqref{e.c.red50} with $\pi_\mu$ replaced by
$\Pi_\mu$. The derivatives of $\Pi_\mu$, after collecting and
introducing 1-forms $\Lambda_I$ containing $\der t_I$,
$I=1,2,3,4$, read
\begin{align}
 \der\Pi_1=& \Lambda_1\w\theta^1+\Pi_8\w\theta^2-\Pi_2\w\theta^4,\nonumber \\
 \der\Pi_2=&  \tfrac{1}{2}\Lambda_1\w\theta^2-\Pi_1\w\Pi_2-\Pi_2\w\Pi_3+\Pi_8\w\theta^3+\frac{u_3^3}{u_1^3}W\Pi_9\w\theta^1\nonumber \\
   &+2f_1\theta^1\w\theta^3+f_4\theta^1\w\theta^4+f_2\theta^2\w\theta^3+f_5\theta^2\w\theta^4,\label{e.c.prol30}\\
 \der\Pi_3=&\tfrac{1}{2}\Lambda_1\w\theta^1+\Pi_8\w\theta^2+\Pi_9\w\theta^3+f_1\theta^1\w\theta^2+f_2\theta^1\w\theta^3+f_5\theta^1\w\theta^4+f_3\theta^2\w\theta^3,\nonumber \\
 \der\Pi_8=&\Lambda_2\w\theta^1+\Lambda_3\w\theta^1+\tfrac{1}{2}\Lambda_1\w\theta^4+\Pi_9\w\Pi_2+\Pi_8\w\Pi_3+f_2\theta^3\w\theta^4,\nonumber\\
 \der\Pi_9=&\Lambda_3\w\theta^1+\Lambda_4\w\theta^2+\Pi_1\w\Pi_9-2\Pi_3\w\Pi_9+\Pi_8\w\theta^4-f_1\theta^1\w\theta^4+f_3\theta^3\w\theta^4. \nonumber
\end{align}
where $f_1,f_2,f_3,f_4,f_5$ are functions on
$G^{prol}\times\P^{(3)}$. We choose the subbundle $\P^c$ of
$G^{prol}\times\P^{(3)}$ by the condition that $f_1,f_2,f_3$ are
equal to zero on $\P^c$. This is done by appropriate specifying of
parameters $t_2,t_3,t_4$ as functions of
$(x$,$y$,$p$,$q$,$u_1,u_2$,$u_3$,$u_8$,$u_9$,$t_1)$. We skip
writing these complicated formulae. The structural equations on
$\P^c$ read \begin{align}
  \der\theta^1 =&\Pi_1\w\theta^1+\theta^4\w\theta^2,\notag \\
  \der\theta^2 =&\Pi_2\w\theta^1+\Pi_3\w\theta^2+\theta^4\w\theta^3, \notag \\
  \der\theta^3 =&\Pi_2\w\theta^2+(2\Pi_3-\Pi_1)\w\theta^3+A\,\theta^4\w\theta^1, \notag \\
  \der\theta^4 =&\Pi_8\w\theta^1+\Pi_9\w\theta^2+(\Pi_1-\Pi_2)\w\theta^4, \notag \\
  \der\Pi_1 =&\Lambda_1\w\theta^1+\Pi_8\w\theta^2-\Pi_2\w\theta^4, \notag \\
  \der\Pi_2 =&(\Pi_3-\Pi_1)\w\Pi_2+A\,\Pi_9\w\theta^1+\tfrac{1}{2}\Lambda_1\w\theta^2
   +\Pi_8\w\theta^3+B\,\theta^1\w\theta^4+C\,\theta^2\w\theta^4, \notag  \\
  \der\Pi_3 =&\tfrac{1}{2}\Lambda_1\w\theta^1+\Pi_8\w\theta^2+\Pi_9\w\theta^3
   +C\,\theta^1\w\theta^4,\label{e.c.prol40} \\
  \der\Pi_8 =&\Pi_9\w\Pi_2+\Pi_8\w\Pi_3-2C\,\Pi_9\w\theta^1+\tfrac{1}{2}\Lambda_1\w\theta^4
   +D\,\theta^1\w\theta^2+2E\,\theta^1\w\theta^3 \notag \\
   &+G\,\theta^1\w\theta^4+H\,\theta^2\w\theta^3+J\,\theta^2\w\theta^4,\notag \\
  \der\Pi_9 =&(\Pi_1-2\Pi_3)\w\Pi_9+\Pi_8\w\theta^4 +E\,\theta^1\w\theta^2+H\,\theta^1\w\theta^3
  +J\,\theta^1\w\theta^4 +L\,\theta^2\w\theta^3, \notag \\
  \der\Lambda_1 =&\Lambda_1\w\Pi_1+2\Pi_8\w\Pi_2+2C\,\Pi_8\w\theta^1-2C\,\Pi_9\w\theta^2
  -A\,\Pi_9\w\theta^4+\wt{M}\,\theta^1\w\theta^2 \notag \\
   &+2(D+AL)\,\theta^1\w\theta^3+\wt{N}\,\theta^1\w\theta^4+2E\,\theta^2\w\theta^3
   +G\,\theta^2\w\theta^4 \notag
\end{align}
with certain functions $A,B,C,D,E,F,G,H,J,L,\wt{M},\wt{N}$ on
$\P^c$.

Above structural equations \emph{uniquely define the only
remaining auxiliary form $\Lambda_1$}. In this manner we
constructed the bundle $\P^c\to\J^2$ and the fixed coframe
associated to the ODEs modulo contact transformations.

\subsubsection{Cartan normal connection from Tanaka's
theory}\label{s.c.normalcon} The above coframe is not fully
satisfactory from the geometric point of view since it does not
transform equivariantly  along the fibres of $\P^c\to\J^2$.

In order to see this we consider the simplest case, related to the
equation $y'''=0$, when all the functions $A,\ldots,\wt{N}$
vanish. Then \eqref{e.c.prol40} become the Maurer-Cartan equations
for the Lie algebra $\o(3,2)\cong\sp(4,\real)$ and $\P^c$ is
locally the Lie group $SP(4,\real)$. The Maurer-Cartan form on
$\P^c$ in the four-dimensional defining representation of
$\sp(4,\real)$ is given by \ben
 \wt{\omega}=\bma \tfrac{1}{2}\Pi_1 & \tfrac{1}{2}\Pi_{2} & -\tfrac{1}{2}\Pi_{8}& -\tfrac{1}{4}\Lambda_1 \\\\
             \hc{4} & \Pi_3-\tfrac{1}{2}\Pi_{1} & -\Pi_{9} & -\tfrac{1}{2}\Pi_{8} \\\\
             \hc{2} & \hc{3} & \tfrac{1}{2}\Pi_{1}-\Pi_{3} & -\tfrac{1}{2}\Pi_{2} \\\\
             2\hc{1} & \hc{2} & -\hc{4} & -\tfrac{1}{2}\Pi_{1}
        \ema.
\een However, this object is not a Cartan connection in a general
case, when $A,\ldots,\wt{N}$ do not vanish, since its curvature
$\wt{K}=\der \wt{\omega}+\wt{\omega}\w\wt{\omega}$ is not
horizontal  with respect to the fibration $\P\to\J^2$, that is the
value of $\wt{K}$ on a vector tangent to a fibre of $\P\to\J^2$ is
not necessarily zero; for instance $\wt{K}^1_{~2}$ contains the
term $A \Pi_9\w\theta^1$.

In order to resolve this problem H. Sato and Y. Yoshikawa
\cite{Sat} found the structural equations for the normal
connection in this problem by means of the Tanaka theory. We
recalculate their result in our notation and give explicit form of
the normal connection, which their paper does not contain.

Let $E^i_{~j}\in\gl(4,\real)$ denotes the matrix whose
$(i,j)$-component is equal to one and other components equal zero.
We introduce the following base in $\sp(4,\real)$
\be\label{e.c.basis_sp}
\begin{aligned}
 &e_1=2E^4_{~1}, & &e_2=E^3_{~1}+E^4_{~2}, & &e_3=E^3_{~2} \\
 &e_4=E^2_{~1}-E^4_{~3}, & &e_5=\tfrac{1}{2}(E^1_{~1}-E^2_{~2}+E^3_{~3}-E^4_{~4}), &
 &e_6=\tfrac{1}{2}(E^1_{~2}-E^3_{~4}),\\
 &e_7=E^2_{~2}-E^3_{~3}, & &e_8=-\tfrac{1}{2}(E^1_{~3}+E^4_{~2}), & &e_9=-E^2_{~3}, \\
 && &e_{10}= -\tfrac{1}{4}E^1_{~4}. &&
\end{aligned}
\ee In this base the form $\wt{\omega}$ is given by
$$\wt{\omega}=\theta^1e_1+\theta^2 e_2+\theta^3 e_3+\theta^4 e_4+\Pi_1e_5+\Pi_2e_6+\Pi_3e_7+\Pi_8e_8+\Pi_9e_9+\Lambda_1e_{10}.$$
The normal Cartan connection $\wh{\omega}^c$ which we seek reads
$$\wh{\omega}^c=\theta^1e_1+\theta^2e_2+\theta^3e_3+\theta^4e_4
+\Omega_1e_5+\Omega_2e_6+\Omega_3e_7+\Omega_4e_8+\Omega_5e_9+\Omega_6e_{10},$$
where the forms $\Omega_\mu$, unknown yet, are given by
\begin{align*}
 &\Omega_1 = \Pi_1+a_i\theta^i, && \Omega_2 = \Pi_2+b_i\theta^i,
 && \Omega_3 =\Pi_3+c_i\theta^i, \nonumber \\
 &\Omega_4 =\Pi_8+f_i\theta^i, && \Omega_5 =\Pi_9+g_i\theta^i,
 && \Omega_6 = \Lambda_1+h_i\theta^i
\end{align*}
and the functions $a_i,b_i,c_i,f_i,g_i,h_i$ are to be found from
the normality conditions of definition \ref{def.Tanakanorm}. The
algebra $\sp(4,\real)$ has the following grading
\be\label{c.gradJ2}
\sp(4,\real)=\g_{-3}\oplus\g_{-2}\oplus\g_{-1}\oplus\g_0\oplus\g_1\oplus\g_2\oplus\g_3,
\ee where \ben
\begin{aligned}
&\g_{-3}=<e_1>, & &\g_{-2}=<e_2>, & &\g_{-1}=<e_3,e_4>, & & \\
&\g_{0}=<e_5,e_7>, & &\g_{1}=<e_6,e_9>, & &\g_{2}=<e_8>, &
&\g_3=<e_{10}>
\end{aligned}
\een and
$\h=\g_0\oplus\g_1\oplus\g_2\oplus\g_3=<e_5,\ldots,e_{10}>$. Let
lower case Latin indices range from $1$ to $4$ and upper case
Latin indices range from $1$ to $10$ throughout this section.

The structure function $\kappa$ of $\wh{\omega}^c$, defined in
\eqref{e.c.kappa}, decomposes into
$$\kappa=\tfrac{1}{2}\kappa^I_{~ij}\,e_I\otimes e^i\w e^j,$$
%\quad i=1,\ldots,10,\quad \alpha,\beta=1,\ldots,4$$
where $(e^I)$ denotes the basis  dual  to $(e_I)$ and
$\kappa^I_{~ij}=\kappa^I_{~[ij]}$ are functions. Condition i) of
definition \ref{def.Tanakanorm} reads
\begin{align*} &\kappa^1_{~23}=0,& &\kappa^1_{~24}=0, &&\kappa^1_{~34}=0, &&\kappa^2_{~34}=0. \end{align*}

\noindent We read structural constants $[e_I,e_J]=c^K_{~IJ} e_K$
for $\sp(4,\real)$, compute the Killing form $B_{IJ}$ and its
inverse $B^{IJ}$. The operator
$\partial^*\colon\Hom(\wedge^2\m,\g)\to\Hom(\m,\g)$ acts as
follows \ben
\partial^*(e_I\otimes e^i\w e^j)=
\left(2\delta^{[i}_{~m}B^{j]K}c^L_{~KI}-\delta^L_{~I}c^{[i}_{~Km}B^{j]K}\right)
e_L\otimes e^m.\een We apply $\partial^*$ to the basis
$(e_I\otimes e^i \w e^j)$ of $\Hom^1(\wedge^2\m,\g)$ and find that
the condition ii) of definition \ref{def.Tanakanorm},
$\partial^*\kappa=0$, is equivalent to vanishing of the following
combinations of $\kappa^I_{~ij}$:
\begin{gather*}
\kappa^1_{~14},\ \kappa^2_{~14},\ \kappa^2_{~24},\
\kappa^3_{~24},\ \kappa^3_{~34},\ \kappa^4_{~23},\
\kappa^1_{~12}+\kappa^2_{~13}+\kappa^2_{~24},\
2\kappa^1_{~12}-\kappa^4_{~24}-\kappa^5_{~34}, \\
\kappa^1_{~12}-\kappa^3_{~23}-\kappa^7_{~34},\
\kappa^1_{~13}+\kappa^2_{~23},\ \kappa^1_{~13}-\kappa^4_{~34},\
\kappa^2_{~12}+\kappa^4_{~14}+\kappa^5_{~24},\\
\kappa^2_{~12}-\kappa^3_{~13}-\kappa^7_{~24},\
\kappa^2_{~12}+\kappa^5_{~24}-\kappa^6_{~34}-\kappa^7_{~24},\
\kappa^2_{~13}+\kappa^3_{~23}+\kappa^5_{~34}-\kappa^7_{~34},\\
\kappa^2_{~23}+\kappa^4_{~34},\
\kappa^3_{~12}-\kappa^5_{~14}+\kappa^6_{~24}+\kappa^7_{~14},\
\kappa^4_{~12}-\kappa^5_{~13}+2\kappa^7_{~13}+\kappa^9_{~24},\\
\kappa^4_{~12}-\kappa^6_{~23}-\kappa^8_{~34}-\kappa^9_{~24},\
\kappa^4_{~13}+\kappa^7_{~23}-\kappa^9_{~34},\
\kappa^4_{~14}+\kappa^6_{~34}+\kappa^7_{~24},\\
\kappa^4_{~24}-\kappa^5_{~34}+2\kappa^7_{~34},\
2\kappa^5_{~12}-2\kappa^8_{~24}-\kappa^{10}_{~34},\
\kappa^5_{~13}+\kappa^6_{~23}-\kappa^8_{~34},\
\kappa^5_{~14}+\kappa^6_{~24},\\
\kappa^5_{~23}-2\kappa^7_{~23}-\kappa^9_{~34},\
2\kappa^6_{~12}+2\kappa^8_{~14}+\kappa^{10}_{~24},\
\kappa^6_{~13}+\kappa^7_{~12}+\kappa^8_{~24}+\kappa^9_{~14}.
\end{gather*}

Next, we calculate the curvature
$\wh{K}^c=\der\wh{\omega}^c+\wh{\omega}^c\w\wh{\omega}^c$, find
the components of the structure function and put them into the
normality conditions. These are only satisfied if all the
functions $a_i,b_i,c_i,e_i,f_i,g_i$ vanish except for
$c_1,h_1,h_2,h_3,h_4$ which are arbitrary and
\begin{align*}
 &a_1=2c_1, && b_1=\tfrac{2}{3}C,&& b_2=c_1, && f_1=\tfrac{2}{3}J, &&
 f_4=c_1,
\end{align*}
where $C$ and $J$ are the functions in \eqref{e.c.prol40}.
Finally, we obtain from the $e_1$-component of the Bianchi
identity
$\der\wh{K}^c=\wh{K}^c\w\wh{\omega}^c-\wh{\omega}^c\w\wh{K}^c$
that
\begin{align*}
 &c_1=0, && h_1=\tfrac{4}{3}G-\tfrac{2}{3}\wt{X}_4(J),&& h_2=\tfrac{2}{3}J, && h_3=0, &&
 h_4=-\tfrac{2}{3}C,
\end{align*}
where $\wh{X}_4$ is the vector field in the frame
$(\wt{X}_1,\wt{X}_2,\wt{X}_3,\wt{X}_4,\wt{X}_5,\wt{X}_6,\wt{X}_7,\wt{X}_8,\wt{X}_9,\wt{X}_{10})$
dual to
$(\theta^1,\theta^1,\theta^1,\theta^1,\Pi_1,\Pi_2,\Pi_3,\Pi_8,\Pi_9,\Lambda_1)$.
The normal connection $\wh{\omega}^c$ has been constructed. The
last thing we must do is renaming the coordinates $u_8\to u_4$,
$u_9\to u_5$ and choosing $u_6$ appropriately, so that formulae
\eqref{e.c.H6} -- \eqref{e.c.Om0} hold. This finishes the proof of
theorem \ref{th.c.1}.

\subsection{Ten-dimensional bundle $\P^c$}\label{s.c.bundle}
\subsubsection{Curvature} We turn to
discussion of consequences of theorem \ref{th.c.1}. The curvature
\ben\wh{K}^c=\der\wh{\omega}^c+\wh{\omega}^c\w\wh{\omega}^c \een
is given by the non-constant terms in the structural equations for
the coframe $\theta^1,\ldots,\Omega_6$.
\begin{align}
 \der\hc{1} =&\vc{1}\w\hc{1}+\hc{4}\w\hc{2},\nonumber \\
 \der\hc{2} =&\vc{2}\w\hc{1}+\vc{3}\w\hc{2}+\hc{4}\w\hc{3},\nonumber \\
 \der\hc{3}=&\vc{2}\w\hc{2}+(2\vc{3}-\vc{1})\w\hc{3}+\inc{A}{2}\hc{2}\w\hc{1}
    +\inc{A}{1}\hc{4}\w\hc{1}, \nonumber \\
 \der\hc{4} =&\vc{4}\w\hc{1}+\vc{5}\w\hc{2}+(\vc{1}-\vc{3})\w\hc{4},\nonumber \\
 \der\vc{1} =&\vc{6}\w\hc{1}+\vc{4}\w\hc{2}-\vc{2}\w\hc{4},\nonumber \\
 \der\vc{2}=&(\vc{3}-\vc{1})\w\vc{2}+\tfrac{1}{2}\vc{6}\w\hc{2}+\vc{4}\w\hc{3}
  +\inc{A}{3}\,\hc{1}\w\hc{2}+\inc{A}{4}\hc{1}\w\hc{4}, \nonumber \\
 \der\vc{3}=&\tfrac{1}{2}\vc{6}\w\hc{1}+\vc{4}\w\hc{2}+\vc{5}\w\hc{3}
  +\inc{A}{5}\hc{1}\w\hc{2}+\inc{A}{2}\hc{1}\w\hc{4},\label{e.c.dtheta_10d} \\
 \der\vc{4}=&\vc{5}\w\vc{2}+\vc{4}\w\vc{3}+\tfrac{1}{2}\vc{6}\w\hc{4}
  +(\inc{A}{6}+\inc{B}{2})\hc{1}\w\hc{2} +2\inc{B}{3}\hc{1}\w\hc{3} \nonumber \\
  &-\inc{A}{3}\hc{1}\w\hc{4}+\inc{B}{4}\hc{2}\w\hc{3} \nonumber \\
 \der\vc{5} =&(\vc{1}-2\vc{3})\w\vc{5}+\vc{4}\w\hc{4}
  +(\inc{A}{7}+\inc{B}{3})\hc{1}\w\hc{2}+\inc{B}{4}\hc{1}\w\hc{3} \nonumber \\
  &-\inc{A}{5}\hc{1}\w\hc{4}+\inc{B}{1}\hc{2}\w\hc{3}, \nonumber \\
  \der\vc{6} =&\vc{6}\w\vc{1}+2\vc{4}\w\vc{2}+\inc{C}{1}\hc{1}\w\hc{2}
   +2\inc{B}{2}\hc{1}\w\hc{3}+\inc{A}{8}\hc{1}\w\hc{4}+2\inc{B}{3}\hc{2}\w\hc{3}, \nonumber
\end{align}
where
$\inc{A}{1},\ldots,\inc{A}{8},\inc{B}{1},\ldots,\inc{B}{4},\inc{C}{1}$
are functions on $\P^c$.

The functions $\inc{A}{1},\ldots,\inc{C}{1}$ are contact relative
invariants of the underlying ODE and the full set of contact
invariants can be constructed by consecutive differentiation of
$\inc{A}{1},\ldots,\inc{C}{1}$ with respect to the frame
$(X_1,X_2,X_3,X_4,X_5,X_6,X_7,X_8,X_9, X_{10})$ dual to
$(\theta^1,$ $\theta^2,$ $\theta^3,$ $\theta^4,$ $\Omega_1,$
$\Omega_2,$ $\Omega_3,$ $\Omega_4,$ $\Omega_5,$ $\Omega_6)$, see
\cite{Olv2}. Utilizing the identities $\der^2\Omega_\mu=0$ we
compute the exterior derivatives of
$\inc{A}{i},\inc{B}{j},\inc{C}{1}$, for instance
$$ \der\inc{B}{1}=X_1(\inc{B}{1})\hc{1}+X_2(\inc{B}{2})\hc{2}+X_3(\inc{B}{3})\hc{3}
-2\inc{B}{4}\hc{4}+2\inc{B}{1}\vc{1}-5\inc{B}{1}\hc{3}.$$ From
these formulae it follows that i) $\inc{A}{2},\ldots\inc{A}{8}$
express by the coframe derivatives of $\inc{A}{1}$, ii)
$\inc{B}{2},\ldots\inc{B}{4}$ express by coframe derivatives of
$\inc{B}{1}$ and iii) $\inc{C}{1}$ is a function of coframe
derivatives of both $\inc{A}{1}$ and $\inc{B}{1}$. Hence we have
\begin{corollary}\label{c.cor.har}
There are two basic contact relative invariants\footnote{This
property means in the language of Tanaka's theory that curvature
of a normal connection is generated by its harmonic part.} for
third order ODEs: \ben \inc{A}{1}=\frac{u_3^3}{u_1^3}W \qquad
\inc{B}{1}=\frac{u_1^2}{6u_3^5}F_{qqqq}. \een All other invariants
can be derived from them by consecutive differentiation with
respect to the dual frame $(X_i)$.
\end{corollary}

The simplest case, in which all the contact invariants
$\inc{A}{1},\ldots,\inc{C}{1}$ vanish corresponds to
$W=F_{qqqq}=0$ and is characterized by
\begin{corollary}\label{cor.c.10d_flat}
For a third-order ODE $y'''=F(x,y,y',y'')$ the following
conditions are equivalent.
\begin{itemize}
 \item[i)] The ODE is contact equivalent to $y'''=0$.
 \item[ii)] It satisfies the conditions $W=0$, and $F_{qqqq}=0$.
  \item[iii)] It has the $\o(3,2)$ algebra of contact symmetry generators.
\end{itemize}
\end{corollary}
%\marginpar{tanaka\_10d.mw}

\subsubsection{Structure of $\P^c$}

The manifold $\P^c$ is endowed with threefold geometry of
principal bundle over the second jet space $\J^2$, the first jet
space $J^1$ and the solution space $\S$. We discuss these
structures consecutively. Let us remind that
$(X_1,X_2,X_3,X_4,X_5,X_6,X_7,X_8,X_9, X_{10})$ denotes the dual
frame to $(\theta^1,$ $\theta^2,$ $\theta^3,$ $\theta^4,$
$\Omega_1,$ $\Omega_2,$ $\Omega_3,$ $\Omega_4,$ $\Omega_5,$
$\Omega_6)$.

First bundle structure, $H_6\to\P\to\J^2$, has been already
introduced explicitly in theorem \ref{th.c.1}. Here we only show
that it is actually defined by the coframe, since it can be
recovered from its structural equations. Indeed, we see from
\eqref{e.c.dtheta_10d} that the ideal spanned by
$\theta^1,\theta^2,\theta^3,\theta^4$ is closed \ben \der
\theta^i\w\theta^1\w\theta^2\w\theta^3\w\theta^4=0 \qquad
\text{for}\qquad i=1,2,3,4, \een and it follows that its
annihilated distribution $<X_5,X_6,X_7,X_8,X_9,X_{10}>$ is
integrable. Maximal integral leaves of this distribution are
locally fibres of the projection $\P^c\to\J^2$. Furthermore, the
commutation relations of these vector fields are isomorphic to the
commutation relations of the six-dimensional group $H_6$, hence we
can {\em define} the fundamental vector fields associated to the
action $H_6$ on $\P^c$ $X_5,\ldots,X_{10}$.

In order to explain how $\P^c$ is the bundle
$CO(2,1)\ltimes\real^3\to\P^c\to\S$ let us first describe the
solution space $\S$ itself. On $\J^2$ there is the congruence of
solutions of the ODE. A family of solutions passing through
sufficiently small open set in $\J^2$ is given by the mapping
$$(x;c_1,c_2,c_3)\mapsto
(x,f(x;c_1,c_2,c_3),f_x(x;c_1,c_2,c_3),f_{xx}(x;c_1,c_2,c_3)),$$
where $y=f(x;c_1,c_2,c_3)$ is the general solution to
$y'''=F(x,y,y',y'')$ and $(c_1,c_2,c_3)$ are constants of
integration. Thus a solution can be considered as a point in the
three-dimensional real space $\S$ parameterized by the constants
of integration. This space can be endowed with a local structure
of differentiable manifold if we {\em choose} a parametrization
$(c_1,c_2,c_3)\mapsto f(x;c_1,c_2,c_3)$ of the solutions and admit
only sufficiently smooth reparameterizations
$(c_1,c_2,c_3)\mapsto(\tilde{c}_1,\tilde{c}_2,\tilde{c}_3)$ of the
constants. We always assume that $\S$ is locally a manifold.
%Hence there is a projection $\J^2\to\S$, whose fibres are solutions
%considered as curves in $\J^2$.
Since $\J^2$ is a bundle over $\S$ so is $\P^c$  and the fibres of
the projection $\P^c\to\S$ are annihilated by the closed ideal
$<\theta^1,\theta^2,\theta^3>$. On the fibres there act the vector
fields $X_4,X_5,X_6,X_7,X_8,X_9,X_{10}$, which form the Lie
algebra $\co(2,1)\semi{.}\real^3$ and thereby define the action of
$CO(2,1)\ltimes\real^3$ on $\P^c$.

Apart from the projection $\J^2\to\S$ there is also the projection
$\J^2\to\J^1$ which takes the second jet $(x,y,p,q)$ of a curve
into its first jet $(x,y,p)$. It gives rise to the third bundle
structure, $H_7\to\P^c\to\J^1$. Here the tangent distribution is
$<X_3,X_5,X_6,X_7,X_8,X_9,X_{10}>$ and it defines the action of a
seven-dimensional group $H_7$ which of course contains $H_6$.

It appears that, under some conditions, $\widehat{\omega}_c$ is
not only a Cartan connection over $\J^2$ but over $\S$ or $\J^1$
also.

\subsection{Conformal geometry on solution space}\label{s.c.conf}

\subsubsection{Normal conformal connection} We remind the notion of
normal conformal connection. Consider $\real^n$ with coordinates
$(x^\mu)$, $\mu=1,\ldots,n$ equipped with the flat metric
$g=g_{\mu\nu}\der x^\mu\otimes\der x^\nu$ of the signature
$(k,l)$, $n=k+l$. The group $Conf(k,l)$ of conformal symmetries of
$g_{\mu\nu}$ consists of
\begin{itemize}
\item[i)] the subgroup $CO(k,l)=\real\times O(k,l)$ containing the
group $O(k,l)$ of isometries of $g$ and the dilatations,
\item[ii)] the subgroup $\real^n$ of translations, \item[iii)] the
subgroup $\real^n$ of special conformal transformations.
\end{itemize}
The stabilizer of the origin in $\real^n$ is the semisimple
product of $CO(k,l)\ltimes\real^n$ of the isometries, the
dilatations, and the special conformal transformations. The flat
conformal space is the homogeneous space
$Conf(k,l)/CO(k,l)\ltimes\real^n$. To this space there is
associated the flat Cartan connection on the bundle
$CO(k,l)\ltimes\real^n\to Conf(k,l)\to \real^n$ with values in the
algebra $\conf(k,l)$.

By virtue of the M\"obius construction, the group $Conf(k,l)$ is
isomorphic to the orthogonal group
$O(k+1,l+1)$ preserving the metric \ben \bma 0 & 0 & -1 \\
0 & g_{\mu\nu} & 0 \\ -1 & 0 & 0 \ema \een on $\real^{n+2}$. This
isomorphism gives rise to the following representation of
$\conf(k,l)\cong\o(k+1,l+1)$ \be\label{e.c.o32} \bma \phi & g_{\nu\rho}\xi^\rho & 0 \\
v^\mu & \lambda^\mu_{~\nu} & \xi^\mu \\ 0 & g_{\nu\rho}v^\rho &
-\phi \ema. \ee Here the vector $(v^\mu)\in\real^n$ generates the
translations, the matrix $(\lambda^\mu_{~\nu})\in\o(k,l)$
generates the isometries, $\phi$ -- the dilatations, and
$(\xi^\mu)\in\real^n$ -- the special conformal transformations.

Let us turn to an arbitrary case of a conformal metric $[g]$ of
the signature $(k,l)$ on a $n$-dimensional manifold $\M$,
$n=k+l>2$. Let us choose a representative $g$ of $[g]$ and
consider a coframe $(\omega^\mu)$, in which
$g=g_{\mu\nu}\,\omega^{\mu}\otimes\omega^\nu$ with constant
coefficients $g_{\mu\nu}$. We calculate the Levi-Civita connection
$\Gamma^\mu_{~\nu}$ for $g$, its Ricci tensor $R_{\mu\nu}$ and the
Ricci scalar $R$. Next we define the following one-forms \ben
P_\nu=\left( \tfrac{1}{2-n}R_{\nu\rho}+\tfrac{1}{2(n-1)(n-2)}R
g_{\nu\rho} \right)\theta^\rho. \een Given these objects we build
the following $\conf(k,l)$-valued matrix $\omega^{conf}$ on $M$
\be\label{e.cnc} \omega^{conf}=\bma
0&P_\nu&0\\
\theta^\mu&\Gamma^\mu_{~\nu}&g^{\mu\rho} P_\rho\\
0&g_{\nu\rho}\theta^\rho&0 \ema. \ee This is the normal conformal
connection on $\M$ in the natural gauge.\footnote{The gauge is
natural since we have started from the Levi-Civita connection, not
from any Weyl connection for $g$, in which case \eqref{e.cnc}
contains a Maxwell potential $A$.} Consider now the conformal
bundle $CO(k,l)\ltimes\real^n\to\P\to\M$, and choose a coordinate
system $(h,x)$ on $\P$ compatible with the local triviality
$\P\cong CO(k,l)\ltimes\real^n \times\M$, where $x$ stands for
$(x^\mu)$ in $\M$ and the matrix $h\in CO(k,l)\ltimes\real^n$
reads \ben h=\begin{pmatrix} {\rm e}^{-\phi}&{\rm
e}^{-\phi}g_{\nu\rho}\xi^\rho &
\frac{1}{2}{\rm e}^{-\phi}\xi^\rho\xi^\sigma g_{\rho\sigma}\\
0&\Lambda^{\mu}_{~\nu}&\Lambda^\mu_{~\rho}\xi^\rho\\
0&0&{\rm e}^\phi
\end{pmatrix},
\quad\quad\Lambda^\mu_{~\rho}\Lambda^\nu_{~\sigma}g_{\mu\nu}=g_{\rho\sigma}.
 \een

The normal conformal connection for $g$ is the following
$\conf(k,l)$-valued one-form on $\P$ \ben
\wh{\omega}^{conf}(h,x)=h^{-1}\cdot\pi^*(\omega^{conf}(x))\cdot h
+h^{-1}\der h. \een The curvature of the normal conformal
connection is as follows \ben
\wh{K}^{conf}(h,x)=h^{-1}\cdot\pi^*(K^{conf}(x))\cdot h, \een
where $K^{conf}$ is the curvature for $\omega^{conf}$ on $\M$ \ben
K^{conf}=\bma
0&D P_\nu&0\\
0& C^\mu_{~\nu} &g^{\mu\rho}D P_\rho\\
0&0&0 \ema, \een and  \ben DP_\mu=\der
P_\mu+P_\nu\w\Gamma^\nu_{~\mu}=\tfrac12P_{\mu\nu\rho}\omega^\nu\w\omega^\rho.
\een The curvature contains the lowest-order conformal invariants
for $g$, namely
\begin{itemize}
\item For $n\geq4$ \ben C^\mu_{~\nu}=\tfrac12
C^\mu_{~\nu\rho\sigma}\omega^\rho\w\omega^\sigma\een is the Weyl
conformal tensor, while \ben
P_{\mu\nu\rho}=\tfrac{1}{n-3}\nabla_{\sigma}C^\sigma_{~\mu\nu\rho}\een
is its divergence. \item For $n=3$ the Weyl tensor identically
vanishes, $C^\mu_{~\nu}= 0$, and the lowest-order conformal
invariant is the Cotton tensor $P_{\mu\nu\rho}$. It has five
independent components.
\end{itemize}

The normality of conformal connections, originally defined by E.
Cartan, is the following property. The algebra
$\conf(k,l)\cong\o(k+1,l+1)$ is graded:
$\o(k+1,l+1)=\g_{-1}\oplus\g_0\oplus\g_1,$ where translations are
the $\g_{-1}$-part, $\co(k,l)$ is the $\g_0$-part and the special
conformal transformations are the $\g_1$-part. The normal
connection for $[g]$ is the only conformal connection such that
the $\co(k,l)$-part of its curvature, $C^\mu_{~\nu}=\tfrac12
C^\mu_{~\nu\rho\sigma}\omega^\rho\w\omega^\sigma$, is traceless:
$C^\rho_{~\nu\rho\sigma}=0$. Cartan normal conformal connections
are normal in the sense of Tanaka.

\subsubsection{Conformal geometries from ODEs}

The following theorem holds.
\begin{theorem}[S.-S.Chern]\label{th.conf3}
If a third-order ODE satisfies a contact-invariant condition
$W=0$, then it has a Lorentzian conformal geometry $[g]$ on its
solution space $S$. Two such geometries constructed from contact
equivalent ODEs are diffeomorphic. In the jet coordinates $[g]$ is
represented by
\begin{align}
g&=(\omega^2)^2-2\omega^1\wt{\omega}^3= \label{e.c.g} \\
&=(\der p-q\der x)^2
 -2(\der y-p\der x)(\der q-F\der x-\tfrac13F_q(\der p-q\der x)+K(\der y-p\der
 x)).\notag\end{align}
The normal conformal connection of this geometry is given by
\be\label{e.c.conn_o}
 \wh{\omega}^c=\bma \vc{3} & -\tfrac{1}{2}\vc{6} & -\vc{4} & -\vc{5} & 0 \\
             \hc{1} & \vc{3}-\vc{1} & -\hc{4} & 0 & -\vc{5} \\
             \hc{2} & -\vc{2} & 0 & -\hc{4} & \vc{4} \\
             \hc{3} & 0 &-\vc{2} & \vc{1}-\vc{3} & -\tfrac{1}{2}\vc{6} \\
             0 & \hc{3} & -\hc{2} & \hc{1} & -\Omega^3
        \ema.
\ee
\begin{proof}
Let us define on $\P^c$ the symmetric two-contravariant tensor
field \ben\wh{g}=(\theta^2)^2-2\theta^1\theta^3 \een of signature
$(++- \,0\,0\,0\,0\,0\,0\,0)$. The degenerate directions of
$\wh{g}$ are precisely those tangent to the fibres of $\P^c\to\S$
\ben \wh{g}(X_i,\cdot)=0,\qquad \text{for}\qquad i=4,5,6,7,8,9,10.
\een The Lie derivatives of $\wh{g}$ along the degenerate
directions are as follows \be\label{e.c.lie1}
L_{X_4}\wh{g}=\frac{u_3^3}{u_1^3}W(\theta^1)^2, \qquad\quad
L_{X_7}\wh{g}=2\wh{g}, \ee and \be\label{e.c.lie2} L_{X_i}\wh{g}=0
\qquad\text{for} \qquad i=5,6,8,9,10. \ee Thus, if only $W=0$, all
the degenerate directions but $X_7$ are isometries of $\wh{g}$
whereas $X_7$ is a conformal symmetry. This allows us to {\em
project} $\wh{g}$ to a Lorentzian conformal metric $[g]$ on the
solution space $\S$. By construction the conformal geometries of
contact equivalent ODEs are equivalent. By construction the
conformal geometries of contact equivalent ODEs are equivalent.

To prove that $\wh{\omega}^c$ is the normal connection we notice
that the condition $W=0$ means $\inc{A}{1}=0$ which causes
$\inc{A}{2}=\ldots=\inc{A}{8}=0$. Thus, structural equations
\eqref{e.c.dtheta_10d} do not contain the non-constant terms
proportional to $\theta^4$ and the curvature $\wh{K}^c$ is
horizontal over $\S$. As a consequence, $\wh{\omega}^c$ is a
connection over $\S$. To assure that it is normal we rearrange
$\wh{\omega}^c$ according to the five-dimensional representation
\eqref{e.c.o32} and we check the normality conditions.
\end{proof}
\end{theorem}

\begin{remark} In order to find the explicit expression for $g$ in a
coordinate system $(c_1,c_2,c_3)$ on $\S$ we would have to find
the general solution $y=f(x;c_1,c_2,c_3)$ of the underlying ODE,
then solve the system $(y=f,\ p=f_{x},\ q=f_{xx})$ with respect to
$c_1,\,c_2,\,c_3$ and substitute these formulae into
\eqref{e.c.g}.
\end{remark}

The conformal grading of $\o(3,2)$ is as follows \be
\o(3,2)=\g_{-1}\oplus\g_0\oplus\g_1, \label{c.gradS} \ee where
\ben \g_{-1}=<e_1,e_2,e_3>, \qquad \g_{0}=<e_4,e_5,e_6,e_7>,
\qquad \g_{1}=<e_8,e_9,e_{10}> \een and the base $(e_i)$ is given
by the decomposition of \eqref{e.c.conn_o}:
\ben\wh{\omega}^c=\theta^1e_1+\theta^2e_2+\theta^3e_3+\theta^4e_4
+\Omega_1e_5+\Omega_2e_6+\Omega_3e_7+\Omega_4e_8+\Omega_5e_9+\Omega_6e_{10}.\een
The curvature is equal to \ben
 \wh{K}^c=\bma 0 & DP_1 & DP_2 & DP_3 & 0 \\
             0 & 0 & 0 & 0 & DP_3 \\
             0 & 0 & 0 & 0 & -DP_2 \\
             0 & 0 & 0 & 0 & DP_1 \\
             0 & 0 & 0 & 0 & 0
        \ema \een
with the following components of the Cotton tensor on $\P^c$
\begin{align*}
 DP_1=&-\tfrac12\inc{C}{1}\,\hc{1}\w\hc{2}-\inc{B}{2}\,\hc{1}\w\hc{3}-\inc{B}{3}\,\hc{2}\w\hc{3}, \notag \\
 DP_2=&-\inc{B}{2}\,\hc{1}\w\hc{2}-2\inc{B}{3}\,\hc{1}\w\hc{3}-\inc{B}{4}\,\hc{2}\w\hc{3},  \\
 DP_3=&-\inc{B}{3}\,\hc{1}\w\hc{2}-\inc{B}{4}\,\hc{1}\w\hc{3}-\inc{B}{1}\,\hc{2}\w\hc{3}. \notag
\end{align*}
Finally, we pull-back these formula to $\J^2$ through $u_1=u_3=1$,
$u_2=u_4=u_5=u_6=0$ and get

\begin{align*}
 DP_1=&(\tfrac12M_p+\tfrac16F_qM_q+\tfrac16F_{qqq}K_y+K_qL_q-\tfrac16 K^2
 F_{qqqq}+\notag \\
 &+\tfrac16K_qF_{qqy}-\tfrac16F_{qqyy}-\tfrac13F_{qqq}K_qK+\tfrac13 F_{qqy}K)\omega^1\w\omega^2
 \notag \\
 &+\tfrac12\left(M_q-K_{qqq}K-2K_{qq}K_q+K_{qqy}\right)\omega^1\w\wt{\omega}^3+\notag \\
 &-\tfrac12 L_{qq}\omega^2\w\wt{\omega}^3, \\
 DP_2=&\tfrac12\left(M_q-K_{qqq}K-2K_{qq}K_q+K_{qqy}\right)\omega^1\w\omega^2+\notag \\
 &-L_{qq}\omega^1\w\wt{\omega}^3+\tfrac12 K_{qqq}\omega^2\w\wt{\omega}^3, \notag \\
 DP_3=&-\tfrac12 L_{qq}\omega^1\w\omega^2+\tfrac12 K_{qqq}\omega^1\w\wt{\omega}^3
 -\tfrac16F_{qqqq}\omega^2\w\wt{\omega}^3. \notag
\end{align*}
The formulae for the conformal connection and curvature (in a
slightly different notation) are given in \cite{New3}.
\begin{example}
The simplest nontrivial equations with vanishing W\"unschmann
invariant and are those with four-dimensional transitive groups of
contact symmetries. Up to contact transformations they as follows:
($\mu\in\real$) \ben
F=\mu\left(\frac{q^2}{1-p^2}-p^2+1\right)^{3/2}-\frac{3q^2p}{1-p^2}+p^3-p^2,\een
\ben F=\mu\frac{(2qy-p^2)^{3/2}}{y^2},  \een \ben
F=4\mu(q-p^2)^{3/2}+6qp-4p^3, \een \ben
F=\mu\left(\frac{q^2}{p^2}+p^2\right)^{3/2}+3\frac{q^2}{p}+p^3,
\een \ben F=(q^2+1)^{3/2}, \qquad F=q^{3/2}. \een First of these
examples has the symmetry group $CO(3)$ and three next have
symmetry $CO(1,2)$.
\end{example}

\subsection{Contact projective geometry on first jet
space}\label{s.c-p} \noindent The connection $\wh{\omega}^c$ gives
rise to not only the above conformal structure but also a geometry
on the first jet space $\J^1$. As we know, see corollary
\ref{c.cor.har}, there are two basic contact invariants in the
curvature $\wh{K}^c$: $W$ and $F_{qqqq}$. The condition $W=0$
yields the conformal geometry. Let us now examine the second
possibility
$$F_{qqqq}=0.$$
The above condition yields
$\inc{B}{1}=\inc{B}{2}=\inc{B}{3}=\inc{B}{4}=0$, which removes all
$\theta^i\w\theta^3$ terms in the curvature and turns
$\wh{\omega}^c$ into a $\sp(4,\real)$ Cartan connection on
$H_7\to\P^c\to\J^1$, since in the curvature there are only terms
proportional to $\theta^1,\theta^2,\theta^4$, horizontal with
respect to $\P^c\to\J^1$. A natural question is to what geometric
structure  $\wh{\omega}^c$ is now related. This geometry is the
contact projective structure on $\J^1$ generated by the family of
solutions of the ODE.

\subsubsection{Contact projective geometry}\label{s.cp3} This
geometry has been exhaustively analyzed in \cite{Fox}, see also
\cite{Cap,Cap2}. We will not discuss the general theory here but
focus on an application of the three-dimensional case to the ODEs.
The definition of contact-projective geometry, see D. Fox
\cite{Fox}, adapted to our situation is the following.

\begin{definition}\label{def.c.cp}
A contact projective structure on the first jet space $\J^1$ is
given by the following data.
\begin{itemize}
\item[i)] The contact distribution $\C$, that is the distribution
annihilated by $$\omega^1=\der y -p\der x.$$ \item[ii)] A family
of unparameterized curves everywhere tangent to $\C$ and such
that: a) for a given point and a direction in $\C$ there is
exactly one curve passing through that point and tangent to that
direction, b) curves of the family are among unparameterized
geodesics for some linear connection on $\J^1$.
  \end{itemize}
\end{definition}

A contact projective structure on $\J^1$ is equivalently given by
a family of linear connections, whose geodesic spray contains the
family of curves. For $\nabla$ to belong to this class one needs
\be\label{e.cp.geod}\nabla_V V=\lambda (V) V \ee along every curve
in the family, where $X$ denotes a tangent field to the considered
curve and $\lambda(X)$ is a function. Given two such connections
$\nabla$ and $\wt{\nabla}$, their difference is a $(2,1)$-type
tensor field $$A(X,Y)=\wt{\nabla}_X Y-\nabla_X Y,$$ for all $X$
and $Y$. Simultaneously we have  $A(V,V)=\mu(V) V$ for $V\in\C$,
where $\mu(V)=\wt{\lambda}(V)-\lambda(V)$ and $\mu_w$ at a point
$w\in\J^1$ is a covector on the vector space $\C_w$. By
polarization we obtain \be\label{e.c.A}
A(X,Y)+A(Y,X)=\mu(X)Y+\mu(Y)X, \qquad X,Y\in\C.\ee The connections
associated to a contact projective structure, when considered at a
point, form an affine space characterized by the above $A$.

Let us describe $A$ more explicitly. We choose a frame
$(e_1=\partial_y, e_2=\partial_p, e_3=\partial_x+p\partial_y)$ and
denote the dual frame by $(\sigma^1,\sigma^2,\sigma^3)$. In
particular $\omega^1=\sigma^1$ and $\C=<e_2,e_3>$. Let
$i,j,\ldots=1,2,3$ and $I,J,\ldots=2,3$. Now $\nabla_j e_i=
\Gamma^k_{~ij} e_k$,
$A^k_{~ij}=\wt{\Gamma}^k_{~ij}-\Gamma^k_{~ij}$, $\mu=\mu^I e_I$
and \eqref{e.c.A} reads \be\label{e.c.Aco} A^k_{~(IJ)}=\mu_{(I}
\delta^{k}_{~J)}. \ee  Relevant components are equal to
$A^1_{~22}=A^1_{~(23)}=A^1_{~33}=0$, $A^2_{~22}=\mu_2,$
$A^2_{~(23)}=\tfrac12\mu_3,$ $A^2_{~33}=0,$ $A^3_{~22}=0,$
$A^3_{~(23)}=\tfrac12\mu_2,$ $A^3_{~33}=\mu_3$ and the rest of
$A^k_{~ij}$ is free. Elementary calculations assure us that the
class of admissible connections is a $20$-dimensional subspace of
$27$-dimensional space of all linear connections on $\J^1$.
Another constraint for the connections is given by
\be\label{e.c.ompr} (\nabla_V
\omega^1)V=\nabla_V(\omega^1(V))=0,\qquad V\in\C. \ee In our frame
this is equivalent to $\Gamma^k_{(IJ)}=0$. Combining eq.
\eqref{e.c.Aco} and \eqref{e.c.ompr} we obtain
\begin{proposition}\label{prop.cp.coord}
The following quantities are invariant with respect to a choice of
a connection in the class distinguished by a contact projective
structure on $\J^1$
\begin{subequations}\label{e.cp.niezm}
\begin{align}
&\Gamma^1_{~22}=0, && \Gamma^1_{~(23)}=0, && \Gamma^1_{~33}=0, && \label{e.cp.niezm1} \\
&\Gamma^3_{~22}, && 2\Gamma^3_{~(23)}-\Gamma^2_{~22}, &&
\Gamma^3_{~33}-2\Gamma^2_{~(23)}, && \Gamma^2_{~33}.
\label{e.cp.niezm2}
\end{align}
\end{subequations}
The connection coefficients are calculated in a frame $(e_i)$ such
that $\C=<e_2,e_3>$.

Values of four the unspecified combinations \eqref{e.cp.niezm2}
define a contact projective structure.
\end{proposition}

Among the above connections there is a distinguished subclass of
those connections which covariantly preserve the distribution
$\C$. We shall call them compatible connections. They satisfy not
only \eqref{e.c.ompr} but a stronger condition \ben \nabla_X
\omega^1 =\phi(X) \omega^1, \qquad \text{for all } X,\een with
some one-form $\phi$. Since $\omega^1$ is non-closed a compatible
connection has nonvanishing torsion.

\subsubsection{Contact projective geometries from ODEs} It is
obvious that the family of solutions of an arbitrary third-order
ODE satisfies the conditions i) and ii a) of definition
\ref{def.c.cp}. (Condition ii a) is satisfied with the possible
exception of the direction $\partial_p$, which belongs to $\C$ but
it is not tangent to any solution in general. However, this
exception is irrelevant since our consideration is local on
$T\J^1$.) We ask when the solutions form a subfamily of geodesics
for a linear connection.
\begin{lemma} \label{lem.c-p}
A third-order ODE $y'''=F(x,y,y',y'')$ defines a
contact-projective structure on $\J^1$ if and only if
$F_{qqqq}=0$. Moreover, the quantities \eqref{e.cp.niezm2} are
given by \be\label{e.cp.gam}\bal &\Gamma^3_{~22}=a_3, &&
2\Gamma^3_{~(23)}-\Gamma^2_{~22}=a_2, \\
&\Gamma^3_{~33}-2\Gamma^2_{~(23)}=a_1, && \Gamma^2_{~33}=-a_0,
\eal\ee where \ben F=a_3q^3+a_2q^2+a_1q+a_0. \een
\begin{proof}
The field $V=\tfrac{\der}{\der x}$ tangent to a solution
$(x,f(x),f'(x))$ equals $V=f''e_2+e_3$ in the frame $(e_i)$. The
geodesic equations \eqref{e.cp.geod} read
\begin{align*}
 & (f'')^2\,\Gamma^1_{~22}+2f''\Gamma^1_{~(23)}+\Gamma^1_{~33}=0, \\
 & f'''+(f'')^2\,\Gamma^2_{~22}+2f''\Gamma^2_{~(23)}+\Gamma^2_{~33}=\lambda(V) f'', \\
 & (f'')^2\,\Gamma^3_{~22}+2f''\Gamma^3_{~(23)}+\Gamma^3_{~33}=\lambda(V).
 \end{align*}
First of these equations is equivalent to \eqref{e.cp.niezm1}.
From the remaining equations we have that
$$
f'''=\Gamma^3_{~22}
{f''}^3+(2\Gamma^3_{~(23)}-\Gamma^2_{~22}){f''}^2
+(\Gamma^3_{~33}-2\Gamma^2_{~(23)})f''-\Gamma^2_{~33},
$$
is satisfied along every solution.
\end{proof}
\end{lemma}

The algebra $\sp(4,\real)\cong\o(3,2)$ has the following grading
(apart from those of \eqref{c.gradJ2} and \eqref{c.gradS})
\begin{align}
&\sp(4,\real)=\g_{-2}\oplus\g_{-1}\oplus\g_0\oplus\g_1\oplus\g_2,
\label{c.gradJ1} \\ \intertext{which reads in the base
\eqref{e.c.basis_sp}:}
&\g_{-2}=<e_1>,  \qquad \g_{-1}=<e_2,e_4 >, \notag \\
&\g_{0}=<e_3,e_5,e_7,e_9>, \notag  \\
&\g_{1}=<e_6,e_8>, \qquad \g_{2}=<e_{10}> \notag.
\end{align}
After calculating Tanaka's normality conditions by the method of
section \ref{s.c.normalcon}, we observe that $\wh{\omega}^c$ is
the normal connection with respect to the grading
\eqref{c.gradJ1}. In this manner we have re-proved a known fact
that to a three-dimensional contact projective geometry there is
associated the unique normal $\sp(4,\real)$-valued Cartan
connection.
\begin{proposition}\label{prop.c-p}
If the contact projective geometry on $\J^1$ exists, then
$\wh{\omega}^c$ of theorem \ref{th.c.1} is the normal Cartan
connection for this geometry.
\end{proposition}

From $\wh{\omega}^c$ one may reconstruct the compatible
connections. To do this we just observe that first, second and
fourth equation of
\eqref{e.c.dtheta_10d} can be written as \ben \bma \der \hp{1} \\
\der \hp{2}\\ \der \hp{4} \ema
+\underbrace{\bma -\Omega_1 & 0 & 0 \\
      -\Omega_2 & -\Omega_3 & \theta^3 \\
      -\Omega_4 & -\Omega_5 & \Omega_3-\Omega_1 \ema}_{\displaystyle
      \wh{\Gamma}}
      \wedge
\bma  \hp{1} \\  \hp{2}\\  \hp{4} \ema = \bma  \hp{4}\w\hp{2}
\\\hp{4}\w\hp{3} \\ 0 \ema. \een
 The three by three matrix denoted by $\wh{\Gamma}$ is the
$\g_0\oplus\g_1$-part of $\wh{\omega}^c$. The following
proposition holds.
\begin{proposition}For any section $s\colon\J^1\to\P^c$ the pull-back
$s^*\wh{\Gamma}$ written in the coframe
$(s^*\theta^1,s^*\theta^2,s^*\theta^4)$ is a connection compatible
with the contact projective geometry.
\begin{proof}
First we choose the section $s_0\colon\J^1\to\P^c$ given by $q=0$,
$u_1=1$, $u_3=1$ and $u_2=u_4=u_5=u_6=0$. We denote
$\Gamma=s^*_0\wh{\Gamma}$. In the coframe
$\sigma^1=s^*_0\theta^1$, $\sigma^2=s^*_0\theta^2$ and
$\sigma^3=s^*_0\theta^4$ we have
$-s^*_0\Omega_3=\Gamma^2_{~2}=\Gamma^2_{~2k}\sigma^k$,
$s^*_0\theta^3=\Gamma^2_{~3}=\Gamma^2_{~3k}\sigma^k$ and so on.
Equations \eqref{e.cp.gam} follow from \eqref{e.c.om} and
\eqref{e.c.Om0}, provided that $F_{qqqq}=0$.

Next we consider an arbitrary section $s\colon\J^1\to\P^c$. In the
local trivialization $\P^c\cong H_7\times \J^1$ we have $\P^c\ni
w=(v,x)$, where $v\in H_7$, $x\in\J^1$ and $s$ is given by
$x\mapsto (v(x),x)$. Now
$s^*\wh{\omega}^c(x)=v^{-1}(x)s^*_0\wh{\omega}^c(x)v(x)
+v^{-1}(x)\der v(x)$, and $s^*\wh{\Gamma}$ is the $\g_0\oplus\g_1$
part of $s^*\wh{\omega}^c$. Since the Lie algebra of $H_7$ is
$\g_0\oplus\g_1\oplus\g_2$, every $v(x)$ in the connected
component of the identity may be written in the form
$v(x)=v_2(x)v_1(x)=\exp{(t_2(x)A_2(x))}\exp{(t_1(x)A_1(x))}$ with
$A_2(x)\in\g_2$ and $A_1(x)\in\g_0\oplus\g_1$. It follows that
\ben
s^*\wh{\omega}^c(x)=v^{-1}(x)_2\left\{v^{-1}_1(x)s^*_0\wh{\omega}^c(x)v_1(x)
+v^{-1}_1(x)\der v_1(x)\right\}v_2(x) +v^{-1}_2(x)\der v_2(x)\een
 But the
$\g_0\oplus\g_1$ part of the quantity in the curly brackets is the
connection $\Gamma=s^*_0\wh{\Gamma}$ written in the coframe
$(s^*\theta^1,s^*\theta^2,s^*\theta^4)$ and $\ad v^{-1}(x)$
transforms it into other compatible connection, according to
\eqref{e.c.Aco}.
\end{proof}
\end{proposition}

\subsection{Six-dimensional conformal geometry in the split
signature}\label{s.c6d} \noindent Until now we have not proposed
any geometric structure, apart from $\wh{\omega}^c$, that could be
associated with an ODE of generic type. Motivated by S.-S. Chern's
construction we would like to build some kind of conformal
geometry starting from an arbitrary ODE, which does not
necessarily satisfy the W\"unschmann condition. We following
theorem holds.
\begin{theorem}
A third order ODE defines a six-dimensional manifold $\M^6$ as
space of integral leaves of distribution
$<X_8,X_9,X_{10},X_5+X_7>$. The manifold $\M^6$ is equipped with a
split signature conformal geometry, whose normal conformal
connection has special holonomy $\sp(4,\real)\semi{.}\real^5$. The
$\sp(4,\real)$ part of this connection is given by
$\wh{\omega}^c$.
\end{theorem}

The rest of the section is devoted to prove the theorem. Let us
define the `inverse' of the symmetric tensor field $
\wh{g}=2\theta^1\theta^3-(\theta^2)^2$ of section \ref{s.c.conf}
to be $\wh{g}_{inv}=\wh{g}^{ij}X_i\otimes X_j=2X_1X_3-(X_2)^2$. We
take the $\o(2,1)$-part of the connection $\wh{\omega}^c$ \ben
\Gamma= \bma
\vc{3}-\vc{1} & -\hc{4} & 0  \\
-\vc{2} & 0 & -\hc{4} \\
0 &-\vc{2} & \vc{1}-\vc{3} \ema, \een and the Levi-Civita symbol
$\epsilon_{ijk}$ in three dimensions. Next we define a new
bilinear form $\wh{\gg}$ on $\P^c$ \ben \wh{\gg}
=\epsilon_{ijk}\,\wh{g}^{kl}\,\theta^i\,\Gamma^j_{~l}. \een The
above method of obtaining  $\wh{\gg}$ of the split degenerate
signature from $\wh{g}$ is called the Sparling procedure
\cite{Nur1}. The new metric reads \be\label{e.c.g33}
\wh{\gg}=2(\Omega_1-\Omega_3)\theta^2-2\Omega_2\theta^1+2\theta^4\theta^3
\ee and was given in \cite{Nur1} in a slightly different context
for the first time. We easily find that its degenerated directions
$X_8,X_9,X_{10},$ and $X_5+X_7$ form an integrable distribution,
so that one can consider the six-dimensional space $\M^6$ of its
integral leaves. The degenerated directions $X_8,X_9,$ and
$X_{10}$ are isometries \ben
L_{X_6}\wh{\gg}=L_{X_8}\wh{\gg}=L_{X_9}\wh{\gg}=0,
 \een
whereas the fourth direction, $X_5+X_7$, is a conformal
transformation  \ben L_{(X_5+X_7)}\wh{\gg}=\wh{\gg}.\een This
allows us to project $\wh{\gg}$ to the split signature conformal
metric $[\gg]$ on $\M^6$ without any assumptions about the
underlying ODE.

It is interesting to study the normal conformal connection
associated to this geometry. Since $\P^c$ is a subbundle of the
conformal bundle over $M^6$, we can calculate the $\o(4,4)$-valued
normal conformal connection \eqref{e.cnc} at once on $\P^c$. It is
as follows.
\begin{equation*}%\label{e.c.conf33}
\wh{\mathbf{w}}=\bma \tfrac12\Omega_1&0&0&\tfrac12\Omega_2&-\tfrac12\Omega_4&-\tfrac12\Omega_6&0&0\\\\  % WIERSZ 1
\Omega_1-\Omega_{{3}}&\Omega_3-\tfrac12\Omega_1&\tfrac12\theta_4&\tfrac12\Omega_2&0&-\mathrm{w}^3_{~5}
&\Omega_5&-\tfrac12\Omega_4\\\\ %WIERSZ 2
-\Omega_2&\Omega_2&\tfrac12\Omega_1&\mathrm{w}^3_{~4}&\mathrm{w}^3_{~5}&0&\Omega_4&-\tfrac12\Omega_6\\\\ %WIERSZ 3
\theta_4&0&0&\Omega_3-\tfrac12\Omega_1&-\Omega_5&-\Omega_4&0&0\\\\ %WIERSZ 4
\theta_2&0&0&\theta_3&\tfrac12\Omega_1-\Omega_3&-\Omega_2&0&0\\\\ %WIERSZ 5
\theta_1&0&0&\tfrac12\theta_2&-\tfrac12\theta_4&-\tfrac12\Omega_1&0&0\\\\ %WIERSZ 6
\theta_3&-\theta_3&-\tfrac12\theta_2&0&-\tfrac12\Omega_2&-\mathrm{w}^3_{~4}&\tfrac12\Omega_1-\Omega_3&
\tfrac12\Omega_{{2}}\\\\%WIERSZ 7
0&\theta_2&\theta_1&\theta_3&\Omega_1-\Omega_3&-\Omega_2&\theta_4&-\tfrac12\Omega_1
\ema,
\end{equation*}
where
\begin{align*}
\mathrm{w}^3_{~4}=&\inc{A}{4}\theta^1+\inc{A}{2}\theta^2+\inc{A}{1}\theta^4,\\
\mathrm{w}^3_{~5}=&\tfrac12\Omega_6+\inc{A}{3}\theta^1+\inc{A}{5}\theta^2+\inc{A}{2}\theta^4.
\end{align*}

It appears that this connection is of very special form. We show
that the algebra of its holonomy group is reduced to
$\sp(4,\real)\semi{.}\real^5$. Let us write down the connection as
\begin{align*}
\wh{\mathbf{w}}=&(\Omega_1-\Omega_3)e_1-\Omega_2e_2+\theta^4e_3+\theta^2e_4+\theta^1e_5+\theta^3e_6+ \\
  &+\Omega_1e_7+\Omega_4e_8+\Omega_5e_9+\Omega_6e_{10}+\mathrm{w}^3_{~5}e_{11}+\mathrm{w}^3_{~4}e_{12},
\end{align*}
where $e_1,\ldots,e_{12}$ are appropriate matrices in $\o(4,4)$.
The space \ben V=\,<e_1,\ldots,e_{12}>\,\subset\o(4,4) \een is not
closed under the commutation relations, however, if we extend $V$
so that it contains three commutators $e_{13}=[e_3,e_{12}]$,
$e_{14}=[e_5,e_{10}]$ and $e_{15}=[e_5,e_{12}]$ then
$<e_1,\ldots,e_{15}>$ is a Lie algebra, a certain semidirect sum
of $\sp(4,\real)$ and $\real^5$. Bases of the factors are the
following: \ben \real^5=<e_1+2e_7-2e_{14}, e_{11}, e_{12}, e_{13},
e_{15}>, \een \ben \sp(4,\real)=<e_2+e_{13}, e_3, e_4, e_5,
e_6-e_{15},e_7,e_8,e_9,e_{10},e_{14}>. \een The matrix of
$\wh{\mathbf{w}}$ can be transformed to the following conjugated
representation, which reveals its structure well \ben \bma
\tfrac{1}{2}\vc{1} & \tfrac{1}{2}\vc{2} & -\tfrac{1}{2}\vc{4} &
-\tfrac{1}{4}\vc{6} &
 2\vc{2} & -2\mathrm{w}^3_{~4} & 2\mathrm{w}^3_{~5} & 0\\\\  % WIERSZ 1
 \hc{4} & \vc{3}-\tfrac{1}{2}\vc{1} & -\vc{5} &-\tfrac{1}{2}\vc{4} &
 4\vc{3}-4\vc{1} & -2\vc{2} & 0 & -2\mathrm{w}^3_{~5} \\\\ %WIERSZ 2
\hc{2} & \hc{3} & \tfrac{1}{2}\vc{1}-\vc{3} & -\tfrac{1}{2}\vc{2}
&
 4\hc{3} & 0 & 2\vc{2} & 2\mathrm{w}^3_{~4} \\\\ %WIERSZ 3
2\hc{1} & \hc{2} & -\hc{4} & -\tfrac{1}{2}\vc{1} &
 0 & -4\hc{3} & 4\vc{1}-4\vc{3} & -2\vc{2}\\\\ %WIERSZ 4
0 & 0 & 0 & 0 & \tfrac{1}{2}\vc{1} & \tfrac{1}{2}\vc{2} & \tfrac{1}{2}\vc{4} & \tfrac{1}{4}\vc{6} \\\\ %WIERSZ 5
0 & 0 & 0 & 0 & \hc{4} & \vc{3}-\tfrac{1}{2}\vc{1} & \vc{5} &\tfrac{1}{2}\vc{4} \\\\ %WIERSZ 6
0 & 0 & 0 & 0 & -\hc{2} & -\hc{3} & \tfrac{1}{2}\vc{1}-\vc{3} & -\tfrac{1}{2}\vc{2}\\\\%WIERSZ 7
0 & 0 & 0 & 0 & -2\hc{1} & -\hc{2} & -\hc{4} & -\tfrac{1}{2}\vc{1}
\ema. \een $\wh{\mathbf{w}}$ has the following block structure in
this representation
\ben \wh{\mathbf{w}}=\bma \wh{\omega}^c & \wh{\tau} \\\\
        0 & -\sigma \wh{\omega}^c \sigma\ema,
\een where \ben
\sigma=\bma 0 & 0 & 0 & 1 \\
               0 & 0 & 1 & 0 \\
           0 & 1 & 0 & 0 \\
           1 & 0 & 0 & 0
       \ema.\een

Surprisingly enough the $\sp(4,\real)$-part of $\wh{\mathbf{w}}$,
given by the diagonal blocks, is totally determined by the
$\sp(4,\real)$ connection $\wh{\omega}^c$. In particular, this
relation holds when $W=0$ and $\wh{\omega}^c$ is a conformal
connection itself. In this case we have rather an unexpected link
between conformal connections in dimensions three and six.
\subsection{Further reduction and geometry on five-dimensional
bundle}\label{s.c.furred} Theorem \ref{th.c.1} is a starting point
for further reduction of the structural group since one can use
the non-constant invariants in \eqref{e.c.dtheta_10d} to eliminate
more variables $u_\mu$. From this point of view third-order ODEs
fall into three main classes:
\begin{itemize}
\item[i)] $W=0$, $F_{qqqq}=0$. This class contains the equations
equivalent to $y'''=0$ and is fully characterized by the corollary
\ref{cor.c.10d_flat}. \item[ii)] $W=0$, $F_{qqqq}\neq 0$. It is
not interesting from the geometric point of view since it does not
contain equations with five-dimensional or larger symmetry groups.
One may prove this by doing full group reduction, see
\cite{Godphd}. \item[iii)] $W\neq 0$. This class leads to a Cartan
connection on a five-dimensional bundle and is studied below.
\end{itemize}

Let us assume that $W\neq 0$ and continue reduction by setting
$\inc{A}{1}=1$, $\inc{A}{2}=0$, which gives \ben
 u_1=\sqrt[3]{W}u_3,  \qquad u_5=\frac{1}{3}\frac{W_q}{\sqrt[3]{W^2}}.
\een At this moment the auxiliary variable $u_6$ which was
introduced by the prolongation becomes irrelevant and may be set
equal to zero $$u_6=0.$$ In second step we choose \ben
u_2=\frac{1}{3}Zu_3 \een and finally \ben
u_4=\frac{1}{9}\frac{W_q}{\sqrt[3]{W^2}}M-\frac{1}{3}\sqrt[3]{W}Z_q.
\een The coframe and the underlying bundle $\P^c$ of theorem
\ref{th.c.1} have been reduced to dimension five according to the
following
\begin{theorem}[S.-S. Chern]\label{th.c.2}
A third-order ODE $y'''=F(x,y,y',y'')$ satisfying the contact
invariant condition $W\neq 0$ and considered modulo contact
transformations of variables, uniquely defines a 5-dimensional
bundle $\P^c_5$ over $\J^2$ and an invariant coframe
$(\theta^1,\ldots,\theta^4,\Omega)$ on it. In local coordinates
$(x,y,p,q,u)$ this coframe is given by \be
\label{e.c.theta_5d}\begin{aligned}
 \theta^1=&\sqrt[3]{W}u\omega^1, \\
 \theta^2=&\frac{1}{3}Z u\omega^1+u\omega^2, \\
\theta^3=&\frac{u}{\sqrt[3]{W}}\left(K+\frac{1}{18}Z^2\right)\omega^1
+\frac{u}{3\sqrt[3]{W}}\left(Z-F_q\right)\omega^2+\frac{u}{\sqrt[3]{W}}\omega^3, \\
\theta^4=&\left(\frac{1}{9}\frac{W_q}{\sqrt[3]{W^2}}Z-\frac{1}{3}\sqrt[3]{W}Z_q\right)\omega^1
+\frac{W_q}{3\sqrt[3]{W^2}}\omega^2+\sqrt[3]{W}\omega^4, \\
 \Omega=& \left(\left(\frac{1}{9}W_q\D Z-\frac{1}{27}W_qZ^2+\frac{1}{9}W_p Z\right)\frac{1}{W}
 -\frac{1}{3}Z_p-\frac{1}{9}F_qZ_q\right)\omega^1 \\
 &+\left(\frac{W_p}{3W}-\frac{1}{3}Z_q\right)\omega^2+
 \frac{W_q}{3W}\,\omega^3+\frac{1}{3}F_q\,\omega^4
 +\frac{\der u}{u}.
 \end{aligned}\ee
where  $\omega^i$ are defined by the ODE via \eqref{e.omega}. The
exterior derivatives of these forms read
\begin{align}
 d\theta^1=&\Omega\w\theta^1-\theta^2\w\theta^4, \nonumber \\
 d\theta^2=&\Omega\w\theta^2+\inc{a}{}\,\theta^1\w\theta^4-\theta^3\w\theta^4,\nonumber\\
 d\theta^3=&\Omega\w\theta^3+\inc{b}{}\,\theta^1\w\theta^2+\inc{c}{}\,\theta^1\w\theta^3
  -\theta^1\w\theta^4+\inc{e}\,\theta^2\w\theta^3+\inc{a}{}\,\theta^2\w\theta^4, \label{e.c.dtheta_5d} \\
 d\theta^4=&\inc{f}{}\,\theta^1\w\theta^2+\inc{g}{}\,\theta^1\w\theta^3+\inc{h}{}\,\theta^1\w\theta^4
  +\inc{k}{}\,\theta^2\w\theta^3-\inc{e}{}\,\theta^2\w\theta^4, \nonumber  \\
 d\Omega=&\inc{l}{}\,\theta^1\w\theta^2+(\inc{f}{}-\inc{a}{}\inc{k}{})\,\theta^1\w\theta^3
 +\inc{m}{}\,\theta^1\w\theta^4+\inc{g}{}\,\theta^2\w\theta^3+\inc{h}{}\,\theta^2\w\theta^4.\nonumber
 \end{align}
\end{theorem}
The basic invariants for \eqref{e.c.dtheta_5d} (i.e. generating
the full set of invariants by consecutive taking of coframe
derivatives) are
$\inc{a}{},\inc{b}{},\inc{e}{},\inc{h}{},\inc{k}{}$:
\begin{align}
\inc{a}{}=& \frac{1}{\sqrt[3]{W^2}}\left(K+\frac{1}{18}Z^2+\frac{1}{9}ZF_q-\frac{1}{3}\D Z\right), \nonumber \\
\inc{b}{}=&\frac{1}{3u\sqrt[3]{W^2}}\bigg(\frac{1}{27}F_{qq}Z^2+\left(K_q-\frac{1}{3}Z_p-\frac{2}{9}F_qZ_q\right)Z+\nonumber\\
   &+\left(\frac{1}{3}\D Z-2K\right)Z_q+Z_y+F_{qq}K-3K_p-K_qF_q-F_{qy}+W_q\bigg), \nonumber \\
\inc{e}{}=&\frac{1}{u}\left(\frac{1}{3}F_{qq}+\frac{1}{W}\left(\frac{2}{9}W_qZ-\frac{2}{3}W_p-\frac{2}{9}W_qF_q \right)\right), \label{e.c.basfun_5d}\\
\inc{h}{}=&\frac{1}{3u\sqrt[3]{W}}\bigg(\left(\frac{1}{9}W_qZ^2-\frac{1}{3}W_pZ+W_y-\frac{1}{3}W_q\D Z\right)\frac{1}{W}+\nonumber  \\
 &+\D Z_q+\frac{1}{3}F_qZ_q\bigg),\nonumber \\
\inc{k}{}=&\frac{1}{u^2\sqrt[3]{W}}\left(\frac{2W_q^2}{9W}-\frac{W_{qq}}{3}\right).
\nonumber
\end{align}

Our next aim is to obtain a Cartan connection. First of all we
study the most symmetric case to find the Lie algebra of a
connection. We assume that all the functions
$\inc{a}{},\ldots,\inc{m}{}$ are constant. Having applied the
exterior derivative to \eqref{e.c.dtheta_5d} we get that
$\inc{b}{},\ldots,\inc{m}{}$ are equal to zero and $\inc{a}{}$ is
an arbitrary real constant $\mu$. In this case the equations
\eqref{e.c.dtheta_5d} become the Maurer-Cartan equations for the
algebra $\real^2\oplus_\mu\real^3$. Straightforward calculations
show that this case corresponds to a general linear equation with
constant coefficients.
\begin{corollary} \label{cor.c.5d_flat}
A third-order ODE is contact equivalent to
$$y'''=-2\mu y'+y,$$ where $\mu$ is an arbitrary constant if and only if
 it satisfies
\begin{align}
1)\quad &W\neq 0, \nonumber \\
2)\quad &\frac{1}{\sqrt[3]{W^2}}
\left(K+\frac{1}{18}Z^2+\frac{1}{9}ZF_q-\frac{1}{3}\D Z\right)=\mu \label{e400} \\
3)\quad &2W_q^2-3W_{qq}W=0. \nonumber
\end{align}
Such an equation has the five-dimensional algebra
$\real^2\semi{\mu}\real^3$ of infinitesimal contact symmetries.
The equations with different constants $\mu_1$ and $\mu_2$ are
non-equivalent.
\begin{proof}
Assume that $\inc{a}{}=\mu$, $\inc{k}{}=0$. It follows from
$\der^2\theta^i=0$ and $\der^2\Omega=0$, that this assumption
makes other functions in \eqref{e.c.dtheta_5d} vanish. Put
$y'''=-2\mu y'+y$ into the formulae of theorem \ref{th.c.2} and
check that it satisfies $\inc{a}{}=\mu$, $\inc{k}{}=0$. Every
equation satisfying $\inc{a}{}=\mu$, $\inc{k}{}=0$ is contact
equivalent to it by virtue of the Cartan equivalence method.
\end{proof}
\end{corollary}
Now we immediately find a family of Cartan connections
\begin{theorem}\label{cor.c.geom5d}
An ODE which satisfies the condition \ben
\frac{1}{\sqrt[3]{W^2}}\left(K+\frac{1}{18}Z^2+\frac{1}{9}ZF_q-\frac{1}{3}\D
Z\right)=\mu\een has the solution space equipped with the
following $\real^2$-valued linear torsion-free connection
\ben%\label{e.c.conn_PtoS}
 \widehat{\omega}_\mu=\bma
             -\vc{} & -\hc{4} & 0 \\
             \mu\hc{4} & -\vc{} & -\hc{4} \\
             -\hc{4} & \mu\hc{4} & -\vc{}
        \ema.
\een Its curvature reads
 \ben
   \bma
      R^1_{~1} & R^1_{~2} & 0 \\
      -\mu R^1_{~2} & R^1_{~1} & R^1_{~2} \\
      R^1_{~2} & -\mu R^1_{~2} & R^1_{~1}
   \ema
 \een
 with
\begin{align}
 R^1_{~1}=&(\mu\inc{g}{}+\inc{k}{})\hc{1}\w\hc{2}+(\mu \inc{k}{}-\inc{g}{})\hc{1}\w\hc{3}
 -\inc{g}{}\hc{2}\w\hc{3}, \nonumber \\
 R^1_{~2}=&-\inc{f}{}\hc{1}\w\hc{2}-\inc{g}{}\hc{1}\w\hc{3}-\inc{k}{}\hc{2}\w\hc{3}. \nonumber
\end{align}
The connection is flat if and only if the related ODE is contact
equivalent to $y'''=-2\mu y'+y$
\begin{proof}
The condition $\inc{a}{}=\mu$ together with its differential
consequences $\inc{b}{}=\inc{c}{}=\inc{e}{}=\inc{h}{}=\inc{m}{}=0$
and $\inc{l}{}=-\inc{k}{}-\mu \inc{g}{}$ is the necessary and
sufficient condition for the curvature of $\widehat{\omega}_\mu$
to be horizontal.
\end{proof}
\end{theorem}

It follows  that every ODE as above has its space of solution equipped with a geometric structure consisting of
\begin{itemize}
\item[i)] Reduction of $\gl(3,\real)$ to $\real^2$ represented by
\ben
   \bma
      a_1 & a_2 & 0 \\
      -\mu a_2 & a_1 & a_2 \\
      a_2 & -\mu a_2 & a_1
   \ema.
\een \item[ii)] A linear torsion-free connection $\Gamma$ taking
values in this $\real^2$.
\end{itemize}
The structure is an example of a geometry with special holonomy.
The algebra $\real^2$ is spanned by the unit matrix and \ben
 \mathrm{m}(\mu) = \bma
      0 & 1 & 0 \\
      -\mu & 0 & 1 \\
      1 & -\mu  & 0
   \ema,
\een whose action on $\S$ is more complicated. Its eigenvalue
equation
$$ \det(u\mathbf{1}-\mathrm{m}(\mu))=u^3+2\mu u-1
$$ is the characteristic polynomial of the linear ODE  $y'''=-2\mu y'+y$. If
$\mu<\tfrac{3}{4}\sqrt[3]{2}$ the polynomial has three distinct
roots and $\mathrm{m}$ is a generator of non-isotropic dilatations
acting along the eigenspaces. If $\mu=\tfrac{3}{4}\sqrt[3]{2}$ the
characteristic polynomial has two roots, one of them double, for
other $\mu$s there is one eigenvalue. The action is diagonalizable
only in the case of three distinct eigenvalues.

\section{Geometries of ODEs modulo point
transformations of variables}\label{ch.point}

\subsection{Cartan connection on seven-dimensional
bundle}\label{s.p.th}  \noindent Following the scheme of reduction
given in section \ref{ch.contact} we construct Cartan connection
for ODEs modulo point transformations.
\begin{theorem}[E. Cartan]\label{th.p.1}
To every third order ODE $y'''=F(x,y,y',y'')$ there are associated
the following data.
\begin{itemize}
\item[i)] The principal fibre bundle $H_3\to\P^p\to\J^2$, where
$\dim\P^p=7$, and $H_3$ is the three-dimensional group \be
\label{e.p.H3}  H_3=\bma \sqrt{u_1}, &
\frac12\frac{u_2}{\sqrt{u_1}}, & 0 & 0 \\\\
 0 & \tfrac{u_3}{\sqrt{u_1}}, & 0 & 0 \\\\
 0 & 0 & \tfrac{\sqrt{u_1}}{u_3}, &
-\tfrac12\tfrac{u_2}{\sqrt{u_1}\,\u_3} \\\\
  0 & 0 & 0 & \tfrac{1}{\sqrt{u_1}}\ema.
\ee \item[ii)] The coframe
$(\theta^1,\theta^2,\theta^3,\theta^4,\Omega_1,\Omega_2,\Omega_3)$,
which defines the $\co(2,1)\oplus_{.}\real^3$-valued Cartan
connection $\wh{\omega}^p$ on $\P^p$ by \be\label{e.p.conn_7d}
 \wh{\omega}^p=\bma \tfrac{1}{2}\vp{1} & \tfrac{1}{2}\vp{2} & 0 & 0 \\\\
             \hp{4} & \vp{3}-\tfrac{1}{2}\vp{1} & 0 & 0 \\\\
             \hp{2} & \hp{3} & \tfrac{1}{2}\vp{1}-\vp{3} & -\tfrac{1}{2}\vp{2} \\\\
             2\hp{1} & \hp{2} & -\hp{4} & -\tfrac{1}{2}\vp{1}
        \ema.
\ee
\end{itemize}
Two 3rd order ODEs $y'''=F(x,y,y',y'')$ and
$y'''=\bar{F}(x,y,y',y'')$ are locally point equivalent if and
only if their associated Cartan connections are locally
diffeomorphic, that is there exists a bundle diffeomorphism
$\Phi\colon\bar{\P}^p\supset\bar{\O}\to\O\subset\P^p$ such that
$$\Phi^*\wh{\omega}^p=\overline{\wh{\omega}^p}.$$

Let $(x,y,p,q,u_1,u_2,u_3)=(x^i,u_\mu)$ be a locally trivializing
coordinate system in $\P^p$. Then the value of $\wh{\omega}^p$ at
the point $(x^i,u_\mu)$ in $\P^p$ is given by \ben
\wh{\omega}^p(x^i,u_\mu)=u^{-1}\,\omega^p\,u+u^{-1}\der u \een
where $u$ denotes the matrix \eqref{e.p.H3} and \ben\omega^p= \bma
\tfrac{1}{2}\vp{1}^0&\tfrac{1}{2}\vp{2}^0 & 0 & 0 \\\\
\wt{\omega}^4 & \vp{3}^0-\tfrac{1}{2}\vp{1}^0 & 0 & 0 \\\\
\omega^2 & \wt{\omega}^3 & \tfrac{1}{2}\vp{1}^0-\vp{3}^0 & -\tfrac{1}{2}\vp{2}^0 \\\\
2\omega^1 & \omega^2 & -\wt{\omega}^4 & -\tfrac{1}{2}\vp{1}^0
\ema\een is the connection $\wh{\omega}^p$ calculated at the point
$(x^i,u_1=1,u_2=0,u_3=1)$.
 The forms $\omega^1,\omega^2,\wt{\omega}^3,\omega^4$ read
\begin{align*}%\label{e.p.om}
\omega^1=&\der y-p\der x, \notag \\
\omega^2=&\der p-q\der x, \\
\wt{\omega}^3=&\der q-F\der x-\tfrac13F_q(\der p-q\der x)+K(\der y-p\der x),\notag \\
\wt{\omega}^4=&\der x +\tfrac16F_{qq}(\der y-p\der x).\notag
\end{align*}
The forms $\vc{1}^0,\ldots,\vc{6}^0$ read
\begin{align*} %\label{e.p.Om0}
\vc{1}^0=&-(3K_q+\tfrac29F_{qq}F_q+\tfrac23F_{qp})\,\omega^1+\tfrac16F_{qq}\omega^2, \notag \\
\vc{2}^0=&\left(L+\tfrac16F_{qq}K\right)\,\omega^1
 -(2K_q+\tfrac19F_{qq}F_q+\tfrac13F_{qp})\,\omega^2
 +\tfrac16F_{qq}\wt{\omega}^3-K\wt{\omega}^4,  \\
\vc{3}^0=&-(2K_q+\tfrac16F_{qq}F_q+\tfrac13F_{qp})\,\omega^1+\tfrac13F_{qq}\,\omega^2
+\tfrac13F_q\wt{\omega}^4.\notag
\end{align*}
\begin{proof}
We begin with the $G_p$-structure on $\J^2$ \eqref{e.G_p}, which
encodes an ODE up to point transformations. In the usual locally
trivializing coordinate system $(x,y,p,q,u_1,\ldots,u_8)$ on
$G_p\times\J^2$ the fundamental form $\hp{i}$ is given by
\begin{align*}
 \hp{1}=& u_1\omega^1,\\
 \hp{2}=& u_2\omega^2+u_3\omega^3, \\
 \hp{3}=& u_4\omega^1+u_5\omega^2+u_6\omega^3,\\
 \hp{4}=& u_8\omega^1+u_7\omega^4.
 \end{align*}
We repeat the procedure of section \ref{s.c.proof}. We choose a
connection by the minimal torsion requirement and then reduce
$G_p\times\J^2$ using the constant torsion property. We
differentiate $\hp{i}$ and gather the $\hp{j}\w\hp{k}$ terms into
\be\label{e.p.red10}\begin{aligned}
\der\hp{1}&=\vp{1}\w\hp{1}+\frac{u_1}{u_3 u_7}\hp{4}\w\hp{2}, \\
 \der\hp{2}&=\vp{2}\w\hp{1}+\vp{3}\w\hp{2}+\frac{u_3}{u_6u_7}\hp{4}\w\hp{3}, \\
 \der\hp{3}&=\vp{4}\w\hp{1}+\vp{5}\w\hp{2}+\vp{6}\w\hp{3},\\
 \der\hp{4}&=\vp{8}\w\hp{1}+\vp{9}\w\hp{2}+\vp{7}\w\hp{4}
\end{aligned}\ee
with the auxiliary connection forms $\vp{\mu}$ containing the
differentials of $u_\mu$ and terms proportional to $\hp{i}$. Then
we reduce $G_p\times\J^2$ by
\ben%\label{e.p.red20}
  u_6=\frac{u_3^2}{u_1},\quad\quad\quad u_7=\frac{u_1}{u_3}.
\een  Subsequently, we get formulae identical to
\eqref{e.c.red_u5}, \eqref{e.c.red_u4}:
\begin{align*}
 u_5=&\frac{u_3}{u_1}\left(u_2-\frac{1}{3}u_3F_q\right), \\
 u_4=&\frac{u^2_3}{u_1}K+\frac{u_2^2}{2u_1}
\end{align*}
and also
\ben
 u_8=\frac{u_1}{6u_3}F_{qq}.
\een  After these substitutions the structural equations for
$\hp{i}$ are the following
\begin{align}%\label{e.p.red50}
 \der\hp{1} =&\vp{1}\w\hp{1}+\hp{4}\w\hp{2},\nonumber \\
 \der\hp{2} =&\vp{2}\w\hp{1}+\vp{3}\w\hp{2}+\hp{4}\w\hp{3},\nonumber \\
 \der\hp{3} =&\vp{2}\w\hp{2}+(2\vp{3}-\vp{1})\w\hp{3}+\inp{A}{1}\hp{4}\w\hp{1},\nonumber \\
 \der\hp{4} =&(\vp{1}-\vp{3})\w\hp{4}+\inp{B}{1}\hp{2}\w\hp{1}
   +\inp{B}{2}\hp{3}\w\hp{1}, \nonumber
 \end{align}
with some functions $\inp{A}{1},\inp{B}{1},\inp{B}{2}$. But now,
in contrast to the contact case, the forms $\vp{1},\vp{2},\vp{3}$
are defined by the above equations without any ambiguity, thus
there is no need to prolong and we have the rigid coframe on the
seven-dimensional bundle $\P^p\to\J^2$.
\end{proof}
\end{theorem}

\subsubsection{Curvature} Further analysis of the coframe of
theorem \ref{th.p.1} is very similar to what we have done in
section \ref{ch.contact}. The curvature of the connection is given
by nonconstant terms in the formulae of exterior differentials of
the coframe
%\marginpar{diff\_p\_7d\_berlin.mw}
\begin{align}
\der\hp{1} =&\vp{1}\w\hp{1}+\hp{4}\w\hp{2},\nonumber \\
\der\hp{2} =&\vp{2}\w\hp{1}+\vp{3}\w\hp{2}+\hp{4}\w\hp{3},\nonumber \\
\der\hp{3}=&\vp{2}\w\hp{2}+(2\vp{3}-\vp{1})\w\hp{3}+\inp{A}{1}\hp{4}\w\hp{1},\nonumber \\
\der\hp{4} =&(\vp{1}-\vp{3})\w\hp{4}+\inp{B}{1}\hp{2}\w\hp{1}
   +\inp{B}{2}\hp{3}\w\hp{1}, \nonumber\\
\der\vp{1}=&-\vp{2}\w\hp{4}+(\inp{D}{1}+3\inp{B}{3})\hp{1}\w\hp{2}
   +(3\inp{B}{4}-2\inp{B}{1})\hp{1}\w\hp{3} \label{e.p.dtheta_7d}  \\
   &+(2\inp{C}{1}-\inp{A}{2})\hp{1}\w\hp{4}-\inp{B}{2}\hp{2}\w\hp{3}, \nonumber \\
\der\vp{2}=&(\vp{3}-\vp{1})\w\vp{2}+\inp{D}{2}\hp{1}\w\hp{2}+(\inp{D}{1}+\inp{B}{3})\hp{1}\w\hp{3}
 +\inp{A}{3}\hp{1}\w\hp{4} \nonumber \\
  &+(2\inp{B}{4}-\inp{B}{1})\hp{2}\w\hp{3}+\inp{C}{1}\hp{2}\w\hp{4}, \nonumber \\
\der\vp{3}=&(\inp{D}{1}+2\inp{B}{3})\hp{1}\w\hp{2}+2(\inp{B}{4}-\inp{B}{1})\hp{1}\w\hp{3}
  +\inp{C}{1}\hp{1}\w\hp{4}-2\inp{B}{2}\hp{2}\w\hp{3}, \nonumber
\end{align}
where
$\inp{A}{1},\inp{A}{2},\inp{A}{3},\inp{B}{1},\inp{B}{2},\inp{B}{3},
\inp{B}{4},\inp{C}{1},\inp{D}{1},\inp{D}{2}$ are functions on
$\P^p$. We have
\begin{corollary}\label{p.cor.har}
The set of basic point relative invariants for third order ODEs is
as follows:
\begin{align*}
 \inp{A}{1}=&\frac{u_3^3}{u^3_1}W, \nonumber \\
 \inp{B}{1}=&\frac{1}{u_3^2}\left(\frac{1}{18}F_{qqq}F_q+\frac{1}{36}F_{qq}^2+\frac{1}{6}F_{qqp}\right)
 -\frac{u_2}{6u_3^3}F_{qqq}, \\
 \inp{C}{1}=&\frac{u_3}{u_1^2}\left(2F_{qq}K+\frac{2}{3}F_qF_{qp}-2F_{qy}+F_{pp} +2W_q\right). \nonumber
\end{align*}
\end{corollary}
All other invariants can be derived from $\inp{A}{1},$
$\inp{B}{1}$ and $\inp{C}{1}$ by consecutive differentiation with
respect to the frame $(X_1,X_2,X_3,X_4,X_5,X_6,X_7)$ be the frame
dual to $(\hp{1},\hp{2},\hp{3},\hp{4},\vp{1},\vp{2},\vp{3})$.
Among these derived invariants $\inp{B}{2}$ and $\inp{B}{4}$ are
important:
\begin{align*}
  \inp{B}{2}=&\frac{u_1}{6u_3^3}F_{qqq}, \\
  \inp{B}{4}=&\frac{1}{u_3^2}
  \left(K_{qq}+\frac19F_{qqq}F_q+\frac13 F_{qqp}+\frac{1}{12}F_{qq}^2 \right).
\end{align*}
\begin{corollary}\label{cor.p.7d_flat}
For a third-order ODE $y'''=F(x,y,y',y'')$ the following
conditions are equivalent.
\begin{itemize}
 \item[i)] The ODE is point equivalent to $y'''=0$.
 \item[ii)] It satisfies the conditions $W=0$, $F_{qqq}=0$, $F_{qq}^2+6F_{qqp}=0$
 and $$2F_{qq}K+\frac{2}{3}F_qF_{qp}-2F_{qy}+F_{pp}=0.$$
  \item[iii)] It has the $\co(2,1)\semi{.}\real^3$ algebra of infinitesimal point symmetries.
\end{itemize}
\end{corollary}

%\section{Geometries on seven-dimensional bundle}

The manifold $\P^p$, like $\P^c$, is equipped with the threefold
structure of principal bundle over $\J^2$, $\J^1$ and $\S$.
\begin{itemize}
\item $\P^p$ is the bundle $H_3\to\P^p\to\J^2$ with the
fundamental fields $X_5,X_6,X_7$. \item It is the bundle
$CO(2,1)\to\P^p\to\S$ with the fundamental fields
$X_4,X_5,X_6,X_7$ \item It is also the bundle $H_4\to\P^p\to\J^1$
with the fundamental fields $X_3,X_5,X_6,X_7$.
\end{itemize}

\subsection{Einstein-Weyl geometry on space of
solutions}\label{s.ew} \noindent We describe in detail the
Einstein-Weyl geometry on the solution space.

\subsubsection{Weyl geometry}\label{s.ew.def} A Weyl geometry on
$\M^n$ is a pair $(g,\phi)$ such that $g$ is a metric of signature
$(k,l)$, $k+l=n$ and $\phi$ is a one-form and they are given
modulo the following transformations \ben \phi\to\phi+\der\lambda,
\qquad\qquad g\to e^{2\lambda} g.\een In particular $[g]$ is a
conformal geometry. For any Weyl geometry there exists the Weyl
connection; it is the unique torsion-free
connection such that \ben% \label{e.W.10}
    \nabla g=2\phi\otimes g.
\een The Weyl connection takes values in the algebra $\co(k,l)$ of
$[g]$. Let $(\omega^\mu)$ be a coframe such that for some $g$ of
$[g]$ is equal $g=g_{\mu\nu}\omega^\mu\otimes\omega^\nu$ with
constant coefficients $g_{\mu\nu}$. The Weyl connection one-forms
$\Gamma^\mu_{~\nu}$ are uniquely defined by the relations
\begin{eqnarray}
 &\der\omega^\mu+\Gamma^\mu_{~\nu}\w\omega^\nu=0, \notag \\
 &\Gamma_{(\mu\nu)}=-g_{\mu\nu}\,\phi, \quad \text{where} \quad
 \Gamma_{ij}=g_{jk}\Gamma^k_{~j}. \notag
\end{eqnarray}
The curvature tensor $R^\mu_{~\nu\rho\sigma}$, the Ricci tensor
$\Ric_{\mu\nu}$ and the Ricci scalar $\R$ of a Weyl connection are
defined as follows
\begin{align*} &R^\mu_{~\nu}=\der\Gamma^\mu_{~\nu}+\Gamma^\mu_{~\rho}\w\Gamma^{\rho}_{~\nu}=
\tfrac12R^\mu_{~\nu\rho\sigma}\omega^\rho\w\omega^\sigma \\
&\Ric_{\mu\nu}=R^\rho_{~\mu\rho\nu}, \\
&\R=\Ric_{\mu\nu}g^{\mu\nu}.
\end{align*}
The Ricci scalar has the conformal weight $-2$, that is it
transforms as $\R\to e^{-2\lambda}\R$ when $g\to e^{2\lambda}g$.
Apart from these objects there is the Maxwell two-form \ben F=\der
\phi, \een which is proportional to the antisymmetric part of the
Ricci tensor.

Einstein-Weyl structures are by definition those Weyl structures
for which the symmetric trace-free part of the Ricci tensor
vanishes \ben
 \Ric_{(\mu\nu)}-\tfrac{1}{n}R\cdot g_{\mu\nu} =0.
\een

\subsubsection{Einstein-Weyl structures from ODEs} One sees from
the system \eqref{e.p.dtheta_7d}  that the pair $(\wh{g},\vp{3})$,
where \ben\wh{g}=2\hp{1}\hp{3}-(\hp{2})^2\een is Lie transported
along the fibres of $\P^p\to\S$ in the following way \ben
L_{X_4}g= \inp{A}{1}(\theta^1)^2,\qquad L_{X_5}\wh{g}=0,\qquad
L_{X_6}\wh{g}=0,\qquad L_{X_7}\wh{g}=2\wh{g},\een and \ben
L_{X_3}\Omega_3=\tfrac12\inp{C}{1}\hp{1}, \een
\ben\label{e.p.lie5} L_{X_j}\Omega_3=0,\qquad\text{for}\qquad
j=5,6,7. \een

\noindent Due to these properties  $(\wh{g},\Omega_3)$ descends
along $\P^p\to\S$ to the Weyl structure $(g,\phi)$ on the solution
space $\S$ on condition that \ben W=0 \een and \be\label{e.p.cart}
\left(\tfrac{1}{3}\D F_q
-\tfrac{2}{9}F_q^2-F_p\right)F_{qq}+\tfrac{2}{3}F_qF_{qp}-2F_{qy}+F_{pp}=0.
\ee These conditions are equivalent to E. Cartan's original
conditions in \cite{Car2}. The conformal metric of the Weyl
structure $(g,\phi)$ coincides with the conformal metric of the
contact case and is represented by
\begin{align*}
g&=2\omega^1\wt{\omega}^3-(\omega^2)^2= \\
&=2(\der y-p\der x)(\der q -\tfrac13F_q\der p+ K\der y
+(\tfrac13qF_q-pK-F)\der x)-(\der p-q\der x)^2,
\end{align*} while the Weyl potential is given by
\ben\phi=-(2K_q+\tfrac19F_{qq}F_q+\tfrac13F_{qp})(\der y-p\der
x)+\tfrac13F_{qq}(\der p-q\der x) +\tfrac13F_q\der x. \een

The Weyl connection for this geometry, lifted to
$CO(2,1)\to\P^p\to\S$, now the bundle of orthonormal frames, reads
\ben
 \Gamma=\bma
             -\vp{1} & -\hp{4} & 0 \\
             -\vp{2} & -\vp{3} & -\hp{4} \\
              0 &-\vp{2} & \vp{1}-2\vp{3}
        \ema.
\een The curvature  is as follows \be\label{e.ewcur}
 (R^\mu_{~\nu})=\bma
R^1_{~1}-F&R^1_{~2}&0\\
R^2_{~1}&-F&R^1_{~2}\\
0&R^2_{~1}&-R^1_{~1}-F \ema \ee with
\begin{eqnarray}
&F=\der\vp{3}=2\inp{B}{3}\hp{1}\w\hp{2}+(2\inp{B}{4}-2\inp{B}{1})\hp{1}\w\hp{3}
-2\inp{B}{2}\hp{2}\w\hp{3},\nonumber  \\
&R^1_{~1}=-\inp{B}{3}\hp{1}\w\hp{2}-\inp{B}{4}\hp{1}\w\hp{3}-\inp{B}{2}\hp{2}\w\hp{3},\nonumber \\
&R^1_{~2}=\inp{B}{1}\hp{1}\w\hp{2}+\inp{B}{2}\hp{1}\w\hp{3},\nonumber \\
&R^2_{~1}=-\inp{B}{3}\hp{1}\w\hp{3}+(\inp{B}{1}-2\inp{B}{4})\hp{2}\w\hp{3}.
\nonumber
\end{eqnarray}
The Ricci tensor reads
\ben
\Ric=\begin{pmatrix}
0&-3\inp{B}{3}&3\inp{B}{1}-5\inp{B}{4}\\
3\inp{B}{3}&2\inp{B}{4}&3\inp{B}{2}\\
-3\inp{B}{1}+\inp{B}{1}&-3\inp{B}{2}&0
\end{pmatrix}
\een and satisfies the Einstein-Weyl equations \ben
 \Ric_{(ij)}=\tfrac{1}{3}R \cdot {g}_{ij}
\een with the Ricci scalar $R=6\inp{B}{4}$. The components of the
curvature in the orthogonal coframe given by $u_1=1, u_2=0, u_3=1$
are as follows
\begin{align*}
\inp{B}{1}=&\tfrac{1}{18}F_{qqq}F_q+\tfrac16F_{qqp}+\tfrac{1}{36}F_{qq}^2,
\\
\inp{B}{2}=&\tfrac16F_{qqq}, \\
\inp{B}{3}=&\tfrac16F_{qqy}-\tfrac13F_{qq}K_q-\tfrac16F_{qqq}K-\tfrac{1}{18}F_{qq}F_{qp}
-\tfrac{1}{54}F^2_{qq}F_q-L_q,\\
\inp{B}{4}=&K_{qq}+\frac19F_{qqq}F_q+\tfrac13F_{qqp}+\frac{1}{12}F_{qq}^2.
\end{align*}

\begin{example}
The simplest nontrivial equations possessing the Einstein-Weyl
geometry are \begin{align*} F=&\frac{3\,q^2}{2\,p}, &
F=&\frac{3q^2p}{p^2+1},\\ F=&\mu\frac{(2qy-p^2)^{3/2}}{y^2},&
F=&q^{3/2}. \end{align*}
\end{example}

\subsection{Geometry on first jet space}\label{s.p-p} \noindent In
section \ref{s.c-p} we described how certain ODEs modulo contact
transformations generate contact projective geometry on $J^1$. The
fact that point transformations form a subclass within contact
transformations suggests that ODEs modulo point transformations
define some refined version of contact projective geometry.
Indeed, the only object that is preserved by point transformations
but is not preserved by contact transformations is the projection
$\J^1\to \text{\em xy plane}$, whose fibres are generated by
$\partial_p$. This motivate us to propose the following
\begin{definition}\label{def.p-p}
A point projective structure on $\J^1$ is a contact projective
structure, such that integral curves of the field $\partial_p$ are
geodesics of the contact projective structure.
\end{definition}
We immediately get
\begin{lemma}
The field $\partial_p$ is geodesic for the contact projective
geometry generated by an ODE provided that the ODE satisfies
$$F_{qqq}=0.$$
\begin{proof}
In the notation of section \ref{s.cp3} we have $\partial_p=e_2$
and, from proposition \ref{prop.cp.coord}, $\Gamma^1_{~22}=0$.
Thus $\nabla_2e_2=\lambda e_2$ iff $\Gamma^3_{~22}=0$, which is
equivalent to $F_{qqq}=0$ by means of lemma \ref{lem.c-p}.
\end{proof}
\end{lemma}

However, the condition $F_{qqq}=0$ is not sufficient for
$\wh{\omega}^p$ to be a Cartan connection for the point projective
structure and we show that there does not exist any simple way to
construct a Cartan connection on $\P^p\to J^1$. The algebra
$\co(2,1)\semi{.}\real^3\subset\o(3,2)$ inherits the following
grading from \eqref{c.gradJ1}
$$\co(2,1)\semi{.}\real^3=\g_{-2}\oplus\g_{-1}\oplus\g_{0}\oplus\g_{1}$$
but it is not semisimple, so the Tanaka method cannot be
implemented. Moreover, the broadest generalization of this method
-- the Morimoto nilpotent geometry handling non-semisimple groups
also fails in this case. This is because the Morimoto approach
requires the algebra $\g$ to be equal to the prolongation of its
non-positive part algebra $\g_{-2}\oplus\g_{-1}\oplus\g_0$. In our
case  the prolongation of $\g_{-2}\oplus\g_{-1}\oplus\g_0$ is
larger then $\co(2,1)\semi{.}\real^3$ and equals precisely
$\o(3,2)$, which yields a contact projective structure. Thus only
the contact case is solved by the methods of the nilpotent
geometry.

Lacking a general theory we must search a Cartan connection in a
more direct way. Consider than an ODE satisfying $F_{qqq}=0$. It
follows that $\inp{B}{2}=0$ and $\inp{B}{4}=\inp{B}{1}$ in
equations \eqref{e.p.dtheta_7d}. We seek four one-forms
\begin{align*}
 \Xi_1=&\vp{1}+a_1\theta^1+a_2\theta^2+a_3\theta^4, \\
 \Xi_2=&\vp{2}+b_1\theta^1+b_2\theta^2+b_3\theta^4, \\
 \Xi_3=&\vp{3}+c_1\theta^1+c_2\theta^2+c_3\theta^4, \\
 \Xi_4=&\hp{3}+f_1\theta^1+f_2\theta^2+f_3\theta^4,
\end{align*}
with yet unknown functions $a_1,\ldots,f_3$, such that the matrix
\ben
\bma \tfrac{1}{2}\Xi_{1} & \tfrac{1}{2}\Xi_{2} & 0 & 0 \\\\
             \hp{4} & \Xi_{3}-\tfrac{1}{2}\Xi_{1} & 0 & 0 \\\\
             \hp{2} & \Xi_{4} & \tfrac{1}{2}\Xi_{1}-\Xi_{3} & -\tfrac{1}{2}\Xi_{2} \\\\
             2\hp{1} & \hp{2} & -\hp{4} & -\tfrac{1}{2}\Xi_{1}
        \ema
\een is a Cartan connection on $\P^p\to\J^1$. Calculating the
curvature for this connection we obtain that the horizontality
conditions yield
$$\der a_1=X_1(a_1)\hp{1}+X_2(a_1)\hp{2}+\inp{B}{1}\hp{3}+X_4(a_1)\hp{4}-a_1\vp{1}-a_2\vp{2}.$$
Unfortunately, none combination of the structural functions
$\inp{A}{1},\ldots,\inp{D}{2}$ and their coframe derivatives of
first order satisfies this condition. Therefore we are not able to
build a Cartan connection for an arbitrary point projective
structure. Moreover, since $\inp{B}{1}$ now equal to
$\frac{1}{u_3^2}(\frac{1}{36}F_{qq}^2+\frac{1}{6}F_{qqp})$ is a
basic relative invariant, it seems to us unlikely that among the
coframe derivatives of $\inp{A}{1},\ldots,\inp{D}{2}$ of any order
there exists a function satisfying the above condition. If such a
function existed it would mean that among the derivatives there is
a more fundamental function from which $\inp{B}{1}$ can be
obtained by differentiation.

Of course, we do have a Cartan connection for the point projective
geometry provided that in addition to $F_{qqq}=0$ the conditions
$\inp{B}{1}=\inp{D}{1}=0$ are imposed. However, the geometric
interpretation of these conditions is unclear.

\subsection{Six-dimensional Weyl geometry in the split
signature}\label{s.w6d} \noindent The construction of the
six-dimensional split signature conformal geometry given in
section \ref{ch.contact} has also its Weyl counterpart in the
point case. A similar construction was done by P. Nurowski,
\cite{Nur1} but he considered the conformal metric, not the Weyl
geometry. Here, apart from the tensor \ben
\wh{\gg}=2(\Omega_1-\Omega_3)\theta^2-2\Omega_3\theta^1+2\theta^4\theta^3
\een of \eqref{e.c.g33}, we also have the one-form \ben
\frac12\vp{3}.\een The Lie derivatives along the degenerate
direction $X_5+X_7$ of $\wh{\gg}$ are \ben
L_{(X_5+X_7)}\wh{\gg}=\wh{\gg}\qquad \text{and}\qquad
L_{(X_5+X_7)}\vp{3}=0.\een In this manner the pair
$(\wh{\gg},\tfrac12\vp{3})$ generates the six-dimensional
split-signature Weyl geometry $(\gg,\phi)$ on the six-manifold
$\M^6$ being the space of integral curves of $X_5+X_7$. The
associated Weyl connection is $\co(3,3)\semi{.}\real^6$-valued and
has the
following form. %on $P_7$
\ben%\label{e.p.weyl33}
\Gamma^\mu_{~\nu}=\bma
0&\tfrac12\hp{4}&\Gamma^1_{~3}&0&\Gamma^1_{~5}
&\Gamma^1_{6}\\\\ %WIERSZ 1
\tfrac12\vp{2}&\tfrac12\vp{1}-\tfrac12\vp{3}&\Gamma^2_{~3}&-\Gamma^1_{~5}&0&\Gamma^2_{~6}\\\\ %WIERSZ 2
\tfrac12\hp{4}&0&\tfrac12\hp{3}-\tfrac12\hp{1}&-\Gamma^1_{~6}&-\Gamma^2_{~6}&0\\\\ %WIERSZ 3
0&-\tfrac12\hp{1}&\tfrac12\hp{3}&-\vp{3}&-\tfrac12\vp{2}&-\tfrac12\hp{4}\\\\ %WIERSZ 4
\tfrac12\hp{1}&0&\tfrac12\hp{2}&-\tfrac12\hp{4}&-\tfrac12\vp{1}-\tfrac12\vp{3}&0\\\\ %WIERSZ 5
-\tfrac12\hp{3}&-\tfrac12\hp{2}&0&-\Gamma^1_{~3}&-\Gamma^2_{~3}&\tfrac12\vp{1}-\tfrac32\vp{3}%WIERSZ 6
\ema, \een where
\begin{align*}
\Gamma^1_{~3}=&\tfrac12\vp{2}+\tfrac12\inp{A}{2}\hp{1},\\
\Gamma^2_{~3}=&\inp{A}{3}\hp{1}+\tfrac12\inp{A}{2}\hp{2}+\inp{A}{1}\hp{4}, \\
\Gamma^1_{~5}=&\inp{D}{1}\hp{1}+\inp{B}{3}\hp{2}+(\tfrac32\inp{B}{4}-\inp{B}{1})\hp{3}
+(\inp{C}{1}-\tfrac12\inp{A}{2})\hp{4}, \\
\Gamma^1_{~6}=&(\tfrac12\inp{B}{4}-\inp{B}{1})\hp{1}-\inp{B}{2}\hp{2},\\
\Gamma^2_{~6}=&(\inp{B}{3}+\inp{D}{1})\hp{1}+\tfrac12\inp{B}{4}\hp{2}+\inp{B}{2}\hp{4}.
\end{align*}
Contrary to section \ref{s.c6d}, this connection seems to have
full holonomy $CO(3,3)\ltimes\real^6$ and it is not generated by
$\wh{\omega}^p$ in any simple way. It is also never Einstein-Weyl.

\subsection{Three-dimensional Lorentzian geometry on solution
space}\label{s.lor3} \noindent The geometries of sections
\ref{s.ew} to \ref{s.w6d} are counterparts of respective
geometries of the contact case. The point classification, however,
contains another geometry, which is new when compared to the
contact case. This is owing to the fact that the Einstein-Weyl
geometry of section \ref{s.ew} has in general the non-vanishing
Ricci scalar, which is a weighted function and can be fixed to a
constant by an appropriate choice of the conformal gauge. Thereby
the Weyl geometry on $\S$ is reduced to a
Lorentzian metric geometry. %Such a reduction cannot be done in
%pure conformal geometry associated with contact transformations,
%since there is no weighted function.

These properties of the Weyl geometry are reflected at the level
of the ODEs by the fact that the equations \be\label{so22}
y'''=\frac{3}{2}\frac{(y'')^2}{y'} \ee
 and \be\label{so4}
y'''=\frac{3y'(y'')^2}{{y'}^2+1} \ee  are contact equivalent to
the trivial $y'''=0$ by means of corollary \ref{cor.c.10d_flat}
but they are mutually {\em non-equivalent} under point
transformations and possess the $\o(2,2)$ and $\o(4)$ algebra of
point symmetries respectively. Both of them generate the same flat
conformal geometry but their Weyl geometries differ. After
calculating equations \eqref{e.ewcur} we see that the only
non-vanishing component of their curvature is the Ricci scalar,
which is negative for the equation \eqref{so22} and positive for
\eqref{so4}. In this circumstances we do another reduction step in
the Cartan algorithm setting the Ricci scalar equal to $\pm 6$
respectively\footnote{We choose $\pm 6$ here to avoid large
numerical factors.}, which means $\inp{B}{4}=\pm1$, and obtain a
six-dimensional subbundle $\P^p_6$ of $\P^p$. The invariant
coframe $(\hp{1},\hp{2},\hp{3},\hp{4},\vp{1},\vp{2})$ yields the
local structure of $SO(2,2)$ or $SO(4)$ on $\P^p_6$ while the
tensor $\wh{g}=2\hp{1}\hp{3}-(\hp{2})^2$ descends to a metric
rather than a conformal class on $\S$ by means of conditions
\be\label{p.liem} L_{X_5}\wh{g}=0, \qquad L_{X_6}\wh{g}=0.\ee The
obtained metrics are locally diffeomorphic to the metrics on the
symmetric spaces $SO(2,2)/SO(2,1)$ or $SO(4)/SO(3)$.

In order to generalize this construction to a broader class of
equations we assume that the Ricci scalar of the Einstein-Weyl
geometry is non-zero \ben
6K_{qq}+\frac23F_{qqq}F_q+2F_{qqp}+\frac12F_{qq}^2 \neq 0\een and
set
$$ u_3=\sqrt{\left|6K_{qq}+\frac23F_{qqq}F_q+2F_{qqp}+\frac12F_{qq}^2\right|}$$
in the coframe of theorem \ref{th.p.1}. The
tensor $\wh{g}$ on $\P^p_6$ projects to the metric $g$ on $\S$
provided that the conditions \eqref{p.liem} still hold, which is
equivalent to \ben W=0 \quad \text{and} \quad
(\D+\tfrac23F_q)\left(
6K_{qq}+\frac23F_{qqq}F_q+2F_{qqp}+\frac12F_{qq}^2\right)=0. \een
%The first condition is needed for the conformal geometry to exist
%and the latter one seems to be a higher-order counterpart of
%\eqref{e.p.cart}.
The Cartan coframe on $\P^p_6$ is then given by
\begin{align}
 \der\hp{1} =&\vp{1}\w\hp{1}-\hp{2}\w\hp{4},\nonumber \\
 \der\hp{2} =&\vp{2}\w\hp{1}+\inv{p}{1}\hp{2}\w\hp{3}-\hp{3}\w\hp{4},\nonumber \\
 \der\hp{3}=&\vp{2}\w\hp{2}-\vp{1}\w\hp{3}+\inv{p}{2}\hp{2}\w\hp{3},\nonumber \\
 \der\hp{4} =&\vp{1}\w\hp{4}+\inv{p}{3}\hp{1}\w\hp{2}
   +\inv{p}{4}\hp{1}\w\hp{3}+\inv{p}{5}\hp{1}\w\hp{4}
   -\tfrac12\inv{p}{2}\hp{2}\w\hp{4}+\inv{p}{1}\hp{3}\w\hp{4}, \nonumber\\
 \der\vp{1} =&-\vp{2}\w\hp{4}+\inv{p}{2}\vp{2}\w\hp{1}
   +\inv{p}{6}\hp{1}\w\hp{2}+\inv{p}{7}\hp{1}\w\hp{3}+\inv{p}{4}\hp{2}\w\hp{3}
   +\inv{p}{5}\hp{2}\w\hp{4},\notag \\
 \der\vp{2}=&-\vp{1}\w\vp{2}+\inv{p}{1}\vp{2}\w\hp{3}+\inv{p}{8}\hp{1}\w\hp{2}
 +\inv{p}{9}\hp{1}\w\hp{3}+\inv{p}{10}\hp{2}\w\hp{3}
  +\inv{p}{5}\hp{3}\w\hp{4}, \nonumber
\end{align}
with some functions $\inv{p}{1},\ldots,\inv{p}{10} $ and the
Levi-Civita connection is given by \ben
 \bma
 \Gamma^1_{~1} & \Gamma^1_{~2} & 0 \\
 \Gamma^2_{~1} & 0 & \Gamma^1_{~2} \\
 0 & \Gamma^2_{~1} & -\Gamma^1_{~1}
 \ema,
\een
where
\begin{align*}
\Gamma^1_{~1}=&-\Omega_1+\tfrac12\inv{p}{2}\theta^2,\\
\Gamma^1_{~2}=&\tfrac12\inv{p}{2}\theta^1-\inv{p}{1}\theta^2-\theta^4,
\\
\Gamma^2_{~1}=&-\Omega_2+\tfrac12\inv{p}{2}\theta^3.
\end{align*}
The curvature reads \ben \bma
 R^1_{~1} & R^1_{~2} & 0 \\
 R^2_{~1} & 0 & R^1_{~2} \\
 0 & R^2_{~1} & -R^1_{~1}
 \ema,
\een

\begin{align*}
R^1_{~1}=&\tfrac12(\inv{p}{9}-\inv{p}{6})\theta^1\w\theta^2+(\tfrac14(\inv{p}{2})^2
-\inv{p}{7})\theta^1\w\theta^3+(\inv{p}{4}+X_2(\inv{p}{1})
+\tfrac12\inv{p}{1}\inv{p}{2})\theta^2\w\theta^3,\\\\
R^1_{~2}=&(\inv{p}{10}-\tfrac12
X_2(\inv{p}{2})-\tfrac14(\inv{p}{2})^2)\theta^1\w\theta^2+
(\inv{p}{4}+X_2(\inv{p}{1})+\tfrac12\inv{p}{1}\inv{p}{2})\theta^1\w\theta^3
\\
&+((\inv{p}{1})^2-X_3(\inv{p}{1}))\theta^2\w\theta^3, \\\\
R^2_{~1}=&-\inv{p}{8}\theta^1\w\theta^2+\tfrac12(\inv{p}{6}-\inv{p}{9})\theta^1\w\theta^3
+(-\inv{p}{10}+\tfrac12X_2(\inv{p}{2})+\tfrac14(\inv{p}{2})^2)\theta^2\w\theta^3.
\end{align*}

\section{Geometry of third-order ODEs modulo fibre-preserving transformations of variables}\label{ch.fp}

\subsection{Cartan connection
on seven-dimensional bundle}\label{s.f}
%\section{Fibre-preserving case}\label{s.f}
\noindent The construction of a Cartan connection for the fibre
preserving case is very similar to its point counterpart. This is
due to the fact that every point symmetry of $y'''=0$ is
necessarily fibre-preserving and, as a consequence, the bundle we
will construct is also of dimension seven. Starting from the
$G_f$-structure of \eqref{e.G_f}, which is given by the forms
\ben\bal
 \hf{1}=&u_1\omega^1, \\
 \hf{2}=&u_2\omega^1+u_3\omega^2,  \\
 \hf{3}=&u_4\omega^1+u_5\omega^2+u_6\omega^3,\\
 \hf{4}=&u_7\omega^4, \\
\eal \een and after the substitutions
\begin{align*}
 u_6=&\frac{u_3^2}{u_1},\quad\quad\quad u_7=\frac{u_1}{u_3}, \\
 u_5=&\frac{u_3}{u_1}\left(u_2-\frac{1}{3}u_3F_q\right), \\
 u_4=&\frac{u^2_3}{u_1}K+\frac{u_2^2}{2u_1}
\end{align*}
we get the following theorem.

\begin{theorem}\label{th.f.1}
To every third order ODE $y'''=F(x,y,y',y'')$ there are associated
the following data.
\begin{itemize}
\item[i)] The principal fibre bundle $H_3\to\P^f\to\J^2$, where
$\dim\P^f=7$, and $H_3$ is the three-dimensional group the same as
in the point case
\ben H_3=\bma \sqrt{u_1}, &\frac12\frac{u_2}{\sqrt{u_1}}, & 0 & 0 \\\\
 0 & \tfrac{u_3}{\sqrt{u_1}}, & 0 & 0 \\\\
 0 & 0 & \tfrac{\sqrt{u_1}}{u_3}, &
-\tfrac12\tfrac{u_2}{\sqrt{u_1}\,\u_3} \\\\
  0 & 0 & 0 & \tfrac{1}{\sqrt{u_1}}\ema,
\een
 \item[ii)] The coframe
$(\hf{1},\hf{2},\hf{3},\hf{4},\vf{1},\vf{2},\vf{3})$, which
defines the $\co(2,1)\oplus_{.}\real^3$-valued Cartan connection
$\wh{\omega}^f$ on $\P^f$ by \be\label{e.f.conn_7d}
 \wh{\omega}^f=\bma \tfrac{1}{2}\vf{1} & \tfrac{1}{2}\vf{2} & 0 & 0 \\\\
             \hf{4} & \vf{3}-\tfrac{1}{2}\vf{1} & 0 & 0 \\\\
             \hf{2} & \hf{3} & \tfrac{1}{2}\vf{1}-\vf{3} & -\tfrac{1}{2}\vf{2} \\\\
             2\hf{1} & \hf{2} & -\hf{4} & -\tfrac{1}{2}\vf{1}
        \ema
\ee \end{itemize} Two 3rd order ODEs $y'''=F(x,y,y',y'')$ and
$y'''=\bar{F}(x,y,y',y'')$ are locally fibre-preserving equivalent
if and only if their associated Cartan connections are locally
diffeomorphic, that is there exists a local bundle diffeomorphism
$\Phi\colon\bar{\P}^f\supset\bar{\O}\to\O\subset\P^f$ such that
$\Phi^*\wh{\omega}^f=\overline{\wh{\omega}^f}$. The value of
$\wh{\omega}^f$ at the point $(x^i,u_\mu)$ in $\P^f$ is given by
\ben \wh{\omega}^f(x^i,u_\mu)=u^{-1}\,{\omega}^f\,u+u^{-1}\der u
\een where $u\in H_3$ and \ben{\omega}^f= \bma
\tfrac{1}{2}\vp{1}^0&\tfrac{1}{2}\vp{2}^0 & 0 & 0 \\\\
\wt{\omega}^4 & \vp{3}^0-\tfrac{1}{2}\vp{1}^0 & 0 & 0 \\\\
\omega^2 & \wt{\omega}^3 & \tfrac{1}{2}\vp{1}^0-\vp{3}^0 & -\tfrac{1}{2}\vp{2}^0 \\\\
2\omega^1 & \omega^2 & -\wt{\omega}^4 & -\tfrac{1}{2}\vp{1}^0
\ema\een is given by
\begin{align*}
\omega^1=&\der y-p\der x, \notag \\
\omega^2=&\der p-q\der x, \notag \\
\wt{\omega}^3=&\der q-F\der x-\tfrac13F_q(\der p-q\der x)+K(\der y-p\der x),\notag \\
\omega^4=&\der x  \\
\vf{1}^0=&-K_q\,\omega^1+\tfrac13F_{qq}\omega^2, \notag \\
\vf{2}^0=&L\,\omega^1-K_q\,\omega^2+\tfrac13F_{qq}\wt{\omega}^3-K\omega^4, \notag \\
\vf{3}^0=&-K_q\,\omega^1+\tfrac13F_{qq}\,\omega^2
+\tfrac13F_q\omega^4.\notag
\end{align*}
\end{theorem}

The exterior differentials of the coframe are equal to
\begin{align}
 \der\hf{1} =&\vf{1}\w\hf{1}+\hf{4}\w\hf{2}+\inf{B}{1}\hf{1}\w\hf{2},\nonumber \\
 \der\hf{2} =&\vf{2}\w\hf{1}+\vf{3}\w\hf{2}+\hf{4}\w\hf{3}+\inf{B}{1}\hf{1}\w\hf{3},\nonumber \\
 \der\hf{3}=&\vf{2}\w\hf{2}+(2\vf{3}-\vf{1})\w\hf{3}+\inf{A}{1}\hf{4}\w\hf{1}+\inf{B}{1}\hf{2}\w\hf{3},\nonumber \\
 \der\hf{4} =&(\vf{1}-\vf{3})\w\hf{4}, \nonumber\\
 \der\vf{1}
 =&-\vf{2}\w\hf{4}+(\inf{D}{1}-\inf{B}{2})\hf{1}\w\hf{2}
   +\inf{B}{3}\hf{1}\w\hf{3}+ \label{e.f.dtheta_7d}  \\
   &+(2\inf{C}{1}-\inf{A}{2})\hf{1}\w\hf{4}+\inf{B}{4}\hf{2}\w\hf{3}+\inf{B}{5}\hf{2}\w\hf{4}, \nonumber \\
 \der\vf{2}=&(\vf{3}-\vf{1})\w\vf{2}+\inf{D}{2}\hf{1}\w\hf{2}+(\inf{D}{1}-2\inf{B}{2})\hf{1}\w\hf{3}
  +\inf{A}{3}\hf{1}\w\hf{4} \nonumber \\
  &+\inf{B}{6}\hf{2}\w\hf{3}+(\inf{C}{1}-\inf{A}{2})\hf{2}\w\hf{4}+\inf{B}{5}\hf{3}\w\hf{4}, \nonumber \\
 \der\vf{3}=&(\inf{D}{1}-\inf{B}{2})\hf{1}\w\hf{2}+\inf{B}{3}\hf{1}\w\hf{3}
  +(\inf{C}{1}-\inf{A}{2})\hf{1}\w\hf{4} \notag \\
  &+\inf{B}{4}\hf{2}\w\hf{3}+\tfrac12\inf{B}{5}\hf{2}\w\hf{4}, \nonumber
\end{align}
where
$\inf{A}{1},\inf{A}{2},\inf{A}{3},\inf{B}{1},\inf{B}{2},\inf{B}{3},\inf{B}{4},\inf{B}{5},\inf{B}{6},\inf{C}{1},\inf{D}{1},\inf{D}{2}$
are  functions on $\P^f$. All these invariants express by the
coframe derivatives of $\inf{A}{1},\inf{B}{1},\inf{C}{1}$, which
read
\begin{align}
 \inf{A}{1}=&\frac{u_3^3}{u^3_1}W, \nonumber \\
 \inf{B}{1}=&\frac{1}{3u_3}F_{qq}, \nonumber \\
 \inf{C}{1}=&\frac{u_2}{u_1^2}\left(\tfrac19F_{qq}F_q+\tfrac13F_{qp}+K_q\right)+ \nonumber \\
  &+\frac{u_3}{u_1^2}\left(\tfrac23F_{qq}K-\tfrac13K_qF_q-K_p-\tfrac23F_{qy}\right).\nonumber
\end{align}
If $\inf{A}{1}=0$ then $\inf{A}{2}=\inf{A}{3}=0$ and if
$\inf{B}{1}=0$ then $\inf{B}{i}=0$ for $i=2,\ldots,6$. In
particular we have \ben
\der\inf{B}{1}=-\inf{B}{2}\hf{1}+(\inf{B}{6}-\inf{B}{3})\hf{2}-\inf{B}{4}\hf{3}-\inf{B}{5}\hf{4}
-\inf{B}{1}\vf{3}.\een The flat case is given by vanishing of
$\inf{A}{1}$, $\inf{B}{1}$ and $\inf{C}{1}$.

\subsection{Fibre-preserving geometry from point
geometry}\label{s.fp}
\subsubsection{Fibre-preserving versus point objects} An immediate
observation about the fibre-preserving objects -- $\P^f$ and
$\wh{\omega}^f$ -- is that they are closely related to their point
counterparts of theorem \ref{th.p.1}. Indeed, there is a unique
diffeomorphism $\rho\colon\P^p\to\P^f$ such that
 \be\label{e.f.pf1}
\overset{p}{\hp{1}}=\rho^*\overset{f}{\hf{1}},\qquad\qquad
\overset{p}{\hp{2}}=\rho^*\overset{f}{\hf{2}},\qquad\qquad
\overset{p}{\hp{3}}=\rho^*\overset{f}{\hf{3}}. \ee It is given by
the identity map in the coordinate systems of theorems
\ref{th.p.1} and \ref{th.f.1}. The remaining one-forms are
transported as follows.\be\label{e.f.pf2} \begin{aligned}
\overset{p}{\hp{4}}=&\rho^*(\overset{f}{\hf{4}}+\tfrac12\inf{B}{1}\,\overset{f}{\hf{1}}), \\
\overset{p}{\vp{1}}=&\rho^*(\overset{f}{\vf{1}}+\inf{B}{5}\,\overset{f}{\hf{1}}
-\tfrac12\inf{B}{1}\,\overset{f}{\hf{2}}),\\
\overset{p}{\vp{2}}=&\rho^*(\overset{f}{\vf{2}}+\tfrac12\inf{B}{5}\,\overset{f}{\hf{2}}
-\tfrac12\inf{B}{1}\,\overset{f}{\hf{3}}), \\
\overset{p}{\vp{3}}=&\rho^*(\overset{f}{\vf{3}}+\tfrac12\inf{B}{5}\,\overset{f}{\hf{1}}),
\end{aligned}\ee
where \ben
\inf{B}{5}=-\frac{1}{3u_1}(\D(F_{qq})+\tfrac13F_{qq}F_q). \een The
above formulae enable us to transform easily the fibre-preserving
coframe into the point coframe.  Given the fibre-preserving
coframe $(\overset{f}{\hf{1}},\ldots,\overset{f}{\vf{3}})$ we
compute $\der \overset{f}{\hf{1}}$ and take the coefficient of the
$\overset{f}{\hf{1}}\w\overset{f}{\hf{2}}$ term. This is the
function $\inf{B}{1}$. Next we compute $\der \inf{B}{1}$,
decompose it in the fibre-preserving coframe, take minus function
that stands next to $\overset{f}{\hf{4}}$ and this is
$\inf{B}{5}$. We substitute these functions together with
$(\overset{f}{\hf{1}},\ldots,\overset{f}{\vf{3}})$ into the right
hand side of \eqref{e.f.pf1} and \eqref{e.f.pf2}, where $\rho$ is
the identity transformation of $\P^f$ and the point coframe is
explicitly constructed on $\P^f$.

Now let us consider the inverse construction, from the
fibre-preserving case to the point case. If we have  only the
point coframe $(\overset{p}{\hp{1}},\ldots,\overset{p}{\vp{3}})$
then we can not utilize eq. \eqref{e.f.pf2} since we are not able
to construct the function $\inf{B}{1}$, which is not a point
invariant\footnote{For example the point transformation
$(x,y)\to(y,x)$ destroys the condition $F_{qq}=0$.} and, as such,
does not appear among functions $\inp{A}{1},\ldots,\inp{D}{2}$ in
\eqref{e.p.dtheta_7d} or among their derivatives. However, if we
consider the point coframe {\em and} the function $\inf{B}{1}$
then the construction is possible, since $\inf{B}{5}$ is given by
the derivative $-X_4(\inf{B}{1})$ along the field $X_4$ of the
{\em point} dual frame. Therefore the passage from the point case
to the fibre-preserving case is possible if we supplement the
connection $\wh{\omega}^p$ with the function $\inf{B}{1}$. This
fact implies that each construction of the point case has its
fibre-preserving counterpart which has an additional object
generated by $\inf{B}{1}$.

\subsubsection{ Counterpart of the Einstein-Weyl geometry on $\S$}
%The counterpart of the three-dimensional Einstein-Weyl geometry
This geometry is constructed in the following way. Let
$(\hf{1},\ldots,\vp{3})$ denotes again the fibre-preserving
coframe. Given the objects $\wh{g}=2\hf{1}\hf{3}-(\hf{2})^2$ and
\ben\wh{\phi}=\vf{3}+\frac12\inf{B}{5}\hf{1}, \een let us also
consider the function $\inf{B}{1}$, and ask under what conditions
the triple $(\wh{g},\wh{\phi},\inf{B}{1})$ can be projected to a
geometry on $\S$. There are two possibilities here, either
$\inf{B}{1}=0$ or $\inf{B}{1}\neq0$. If $\inf{B}{1}=0$ then it is
easy to see that the pair $(\wh{g},\wh{\phi})$ generates the
Einstein-Weyl geometry if only $\inf{A}{1}=\inf{C}{1}=0$, which
means that we are in the trivial case $y'''=0$.

Suppose $\inf{B}{1}\neq0$ then. For the geometry on $\S$ to exist
we need not only the conditions for the Lie transport of $\wh{g}$
and $\wh{\phi}$ but also \be\label{e.lieb}
L_{X_i}\inf{B}{1}=0,\quad \text{for}\quad i=4,5,6,\quad
L_{X_7}\inf{B}{1}=-\inf{B}{1}.\ee If all these conditions are
satisfied then $(\wh{g},\wh{\phi},\inf{B}{1})$ defines on $\S$ the
Einstein-Weyl geometry $(g,\phi)$ of the point case, which is
equipped with an additional object: a weighted function $f$ which
transforms $f\to e^{-\lambda}f$ when $g\to e^{2\lambda}g$ and is
given by the projection of $\inf{B}{1}$. The conditions for
existence of this geometry are $\inf{A}{1}=\inf{B}{5}=0$, that is
\be\label{f.lief} W=0 \qquad\text{and}\qquad
\D(F_{qq})+\tfrac13F_{qq}F_q=0.\ee As usual, the condition $W=0$
guarantees existence of $[g]$ and the other condition yields
\eqref{e.lieb}. The proper Lie transport of $\wh{\phi}$ along
$X_4$ is already guaranteed by the above conditions as their
differential consequence.

\subsubsection{Counterpart of the Weyl geometry on $\M^6$} In the
similar vein we show that the triple
$(\wh{\gg},\tfrac12\wh{\phi},\inf{B}{1})$, where \begin{eqnarray*}
&\wh{\gg}=2(\Omega_1-\Omega_3)\theta^2-2\Omega_3\theta^1+2\theta^4\theta^3,
\\
&\wh{\phi}=\vf{3}+\frac12\inf{B}{5}\hf{1}.
\end{eqnarray*}
 projects to the
six-dimensional split signature Weyl geometry $(\gg,\phi)$ of
section \ref{ch.point} section \ref{s.w6d} equipped with a
function $f$ of conformal weight $-2$.

\subsubsection{Counterpart of Lorentzian geometry on $\S$} Given
$(g,\phi,f)$ on $\S$ it is natural to fix the conformal gauge so
as\footnote{With possible change of the signature to make $f$
positive.} $f=1$. This is equivalent to another substitution \ben
u_3=\tfrac13F_{qq} \een in the Cartan reduction algorithm, which
leads us to the bundle $\P^f_6$ with the following differential
system \be\bal
 \der\hf{1}=&\,\vf{1}\w\hf{1}+\hf{4}\w\hf{2},\\
 \der\hf{2}=&\,\vf{2}\w\hf{1}+\inv{f}{1}\hf{3}\w\hf{2}+\inv{f}{2}\hf{4}\w\hf{2}+\hf{4}\w\hf{3},\\
 \der\hf{3}=&-\vf{1}\w\hf{3}+\vf{2}\w\hf{2}+(2-2\inv{f}{3})\hf{3}\w\hf{2}+\inv{f}{4}\hf{4}\w\hf{1}+2\inv{f}{2}\hf{4}\w\hf{3},\\
 \der\hf{4}=&\,\vf{1}\w\hf{4}+\inv{f}{5}\,\hf{4}\w\hf{1}+(\inv{f}{3}-2)\,\hf{4}\w\hf{2}+\inv{f}{1}\,\hf{4}\w\hf{3},\\
 \der\vf{1}=&\,(2\inv{f}{3}-2)\vf{2}\w\hf{1}-\vf{2}\w\hf{4}+\inv{f}{6}\hf{1}\w\hf{2}+\inv{f}{7}\hf{1}\w\hf{3}+\inv{f}{8}\hf{1}\w\hf{4}
    -\inv{f}{5}\hf{2}\hf{4}, \\
 \der\vf{2}=&\,\vf{2}\w\vf{1}-\inv{f}{1}\vf{2}\w\hf{3}-\inv{f}{2}\vf{2}\w\hf{4}+\inv{f}{9}\hf{1}\w\hf{2}
  +\inv{f}{10}\hf{1}\w\hf{3}+\inv{f}{11}\hf{1}\w\hf{4}+\\
  &+\inv{f}{12}\hf{2}\w\hf{3}+\inv{f}{13}\hf{2}\w\hf{4}-\inv{f}{5}\hf{3}\w\hf{4}. \\
 \eal \ee
 If the conditions \eqref{f.lief}, now equivalent to
 $\inv{f}{4}=0=\inv{f}{2}$, are satisfied then $\wh{g}=2\hf{1}\hf{3}-(\hf{2})^2$
projects to a Lorentzian metric $g$ and \ben
\wh{\phi}=-2\inv{f}{5}\hf{1}+2\inv{f}{3}\hf{2}+2\inv{f}{1}\hf{3}
\een projects to a one-form $\phi$. With the Lorentzian metric
there is associated the Levi-Civita connection
$(\Gamma^\mu_{~\nu})$:
\begin{align*}
\Gamma^1_{~1}=&-\vf{1}+(\inv{f}{3}-1)\hf{2},\\
\Gamma^1_{~2}=&(\inv{f}{3}-1)\hf{1}+\inv{f}{1}\hf{2}-\hf{4},
\\
\Gamma^2_{~1}=&-\vf{2}+(\inv{f}{3}-1)\hf{3}.
\end{align*}
%\begin{align}
%R^1_{~1}=&(\inv{f}{5}-\tfrac{1}{2}\inv{f}{6})\hf{1}\w\hf{2}
%+((\inv{f}{3})^2-4\inv{f}{3}-2\inv{f}{12}+1)\hf{1}\w\hf{3}+\notag\\
%&-(X_2(\inv{f}{1})+\inv{f}{1}\inv{f}{3})\hf{2}\w\hf{3},\nonumber  \\
%R^1_{~2}=&((\inv{f}{3})^2-2\inv{f}{3}+\inv{f}{1}\inv{f}{5}-\inv{f}{12}-1)\hf{1}\w\hf{2}
%-(X_2(\inv{f}{1})+\inv{f}{1}\inv{f}{3})\hf{1}\w\hf{3}+\notag\\
%&-(X_3(\inv{f}{1})+\inv{f}{1}^2)\hf{2}\w\hf{3},\nonumber \\
%R^2_{~1}=&-\inv{f}{9}\hf{1}\w\hf{2}-(\inv{f}{5}-\tfrac{1}{2}\inv{f}{6})\hf{1}\w\hf{3}+\nonumber\\
%&-((\inv{f}{3})^2-2\inv{f}{3}+\inv{f}{1}\inv{f}{5}-\inv{f}{12}-1)\hf{2}\w\hf{3}.\nonumber
%\nonumber
%\end{align}
The covariant derivative of $\phi$ with respect to
$\Gamma^\mu_{~\nu}$ is as follows \ben \phi_{i;j}=\begin{pmatrix}
-\inv{f}{9}-(\inv{f}{5})^2&\tfrac{1}{2}\inv{f}{6}+\inv{f}{5}(\inv{f}{3}-2)&\inv{f}{12}-\inv{f}{3}(\inv{f}{3}-3)\\
\tfrac{1}{2}\inv{f}{6}+\inv{f}{5}\inv{f}{3}&2\inv{f}{12}-2\inv{f}{3}(\inv{f}{3}-2)&X_2(\inv{f}{1})+\inv{f}{1}\\
\inv{f}{12}-\inv{f}{3}(\inv{f}{3}-1)&X_2(\inv{f}{1})-\inv{f}{1}&X_3(\inv{f}{1})
\end{pmatrix}.
\een The one-form $\phi$ and the Ricci tensor satisfy the
following identities
\begin{eqnarray*}
 &\nabla_{(i}\phi_{j)}=-\Ric_{ij}-\phi_i\phi_j+(\phi^k\phi_k+2)g_{ij},\\
 %&\nabla_{[i}\phi_{j]}=-(\star\phi)_{ij},\\
 &\R=2\phi^k\phi_k+6,\\
 &\der\phi=-2\ast\phi.
\end{eqnarray*}
The homogeneous model of this geometry is associated to
$y'''=\tfrac32\tfrac{(y'')^2}{y'}$.

\end{document}